\magnification 1200

%
% ??/epsf.tex (written by Radical Eye Software and copied below)
% defines the macro \epsfbox with one argument,
% the encapsulated PostScript file to include.
% Invoking it causes a \vbox with the natural size of the drawing
% to be inserted at the point of invocation.
% Usually figures are meant to be centered and set off, and possibly
% to have a title and/or a figure number. The macros below do that.
% You should assign values that please you to the variables
% \abovefigskip, \belowfigskip, figtitleskip, figtitlefont,...
% AFTER \inputting the present file: \input /u2/kbi/tex/epsfig
% and BEFORE the first invocation of any of the macros below.
% RESET AT YOUR PLEASURE THE VARIABLES AT THE VERY BOTTOM!
%
% For convenience we make a dimension for figures:
\newdimen\FigSize	\FigSize=.9\hsize % alter at your convenience
%
% For a SCALED HORIZONTALLY CENTERED FIGURE use \epsfig.
% First argument is the horizontal width of the figure, given
% in any way TeX can understand.
% The second argument is an encapsulated PostScript file name (filnam.eps).
% Note the mandatory semicolons between arguments in this example!:
% \epsfig .8\hsize; example.ps;
% will put a centered scaled \vbox of width .8\hsize suitably offset 
% at the point of invocation
\newskip\abovefigskip	\newskip\belowfigskip
\gdef\epsfig#1;#2;{\par\vskip\abovefigskip\penalty -500
   {\everypar={}\epsfxsize=#1\noindent
    \centerline{\epsfbox{#2}}}%
    \vskip\belowfigskip}%
%
% SCALED TITLED EPSFIG HORIZONTALLY CENTERED: \tepsfig.
% First argument is the horizontal width of the figure,
% second an encapsulated PostScript file name,
% third a title for the figure.
% Note the mandatory semicolons between arguments!
% example: \tepsfig5truein; example.ps;{This is a figure}
\newskip\figtitleskip
\gdef\tepsfig#1;#2;#3{\par\vskip\abovefigskip\penalty -500
   {\everypar={}\epsfxsize=#1\noindent
    \vbox
      {\centerline{\epsfbox{#2}}\vskip\figtitleskip
       \centerline{\figtitlefont#3}}}%
    \vskip\belowfigskip}%
%
% SCALED NUMBERED TITLED EPSFIG HORIZONTALLY CENTERED: \nepsfig
% The figure number is automatically increased for every
% invocation of \nepsfig or \nipsfig.
\newcount\FigNr	\global\FigNr=0
\gdef\nepsfig#1;#2;#3{\global\advance\FigNr by 1
   \tepsfig#1;#2;{Figure\space\the\FigNr.\space#3}}%
%
% 
% Often you would rather have TeX decide where to put the figure
% by using \midinsert. Here are macros that do that
% (mnemonics ``ipsfig'' is for ``midInsert PS FIGure'')
% 
% TeX-PLACED SCALED EPSFIG HORIZONTALLY CENTERED: \ipsfig
\gdef\ipsfig#1;#2;{%\goodbreak ??
   \midinsert{\everypar={}\epsfxsize=#1\noindent
	      \centerline{\epsfbox{#2}}}%
   \endinsert}%
%
% TeX-PLACED SCALED TITLED EPSFIG HORIZONTALLY CENTERED: \tipsfig
\gdef\tipsfig#1;#2;#3{\midinsert
   {\everypar={}\epsfxsize=#1\noindent
    \vbox{\centerline{\epsfbox{#2}}%
          \vskip\figtitleskip
          \centerline{\figtitlefont#3}}}\endinsert}%
%
% TeX-PLACED SCALED NUMBERED TITLED EPSFIG HORIZONTALLY CENTERED: \nipsfig
% example: \nipsfigd.9\hsize;example.ps;{This is an example figure}
\gdef\nipsfig#1;#2;#3{\global\advance\FigNr by1%
  \tipsfig#1;#2;{Figure\space\the\FigNr.\space#3}}%
\newread\epsffilein    % file to \read
\newif\ifepsffileok    % continue looking for the bounding box?
\newif\ifepsfbbfound   % success?
\newif\ifepsfverbose   % report what you're making?
\newdimen\epsfxsize    % horizontal size after scaling
\newdimen\epsfysize    % vertical size after scaling
\newdimen\epsftsize    % horizontal size before scaling
\newdimen\epsfrsize    % vertical size before scaling
\newdimen\epsftmp      % register for arithmetic manipulation
\newdimen\pspoints     % conversion factor
\pspoints=1bp          % Adobe points are `big'
\epsfxsize=0pt         % Default value, means `use natural size'
\epsfysize=0pt         % ditto
\def\epsfbox#1{\global\def\epsfllx{72}\global\def\epsflly{72}%
   \global\def\epsfurx{540}\global\def\epsfury{720}%
   \def\lbracket{[}\def\testit{#1}\ifx\testit\lbracket
   \let\next=\epsfgetlitbb\else\let\next=\epsfnormal\fi\next{#1}}%
\def\epsfgetlitbb#1#2 #3 #4 #5]#6{\epsfgrab #2 #3 #4 #5 .\\%
   \epsfsetgraph{#6}}%
\def\epsfnormal#1{\epsfgetbb{#1}\epsfsetgraph{#1}}%
\def\epsfgetbb#1{%
%
%   The first thing we need to do is to open the
%   PostScript file, if possible.
%
\openin\epsffilein=#1
\ifeof\epsffilein\errmessage{I couldn't open #1, will ignore it}\else
%
%   Okay, we got it. Now we'll scan lines until we find one that doesn't
%   start with %. We're looking for the bounding box comment.
%
   {\epsffileoktrue \chardef\other=12
    \def\do##1{\catcode`##1=\other}\dospecials \catcode`\ =10
    \loop
       \read\epsffilein to \epsffileline
       \ifeof\epsffilein\epsffileokfalse\else
%
%   We check to see if the first character is a % sign;
%   if not, we stop reading (unless the line was entirely blank);
%   if so, we look further and stop only if the line begins with
%   `%%BoundingBox:'.
%
          \expandafter\epsfaux\epsffileline:. \\%
       \fi
   \ifepsffileok\repeat
   \ifepsfbbfound\else
    \ifepsfverbose\message{No bounding box comment in #1; using defaults}\fi\fi
   }\closein\epsffilein\fi}%
%
%   Now we have to calculate the scale and offset values to use.
%   First we compute the natural sizes.
%
\def\epsfsetgraph#1{%
   \epsfrsize=\epsfury\pspoints
   \advance\epsfrsize by-\epsflly\pspoints
   \epsftsize=\epsfurx\pspoints
   \advance\epsftsize by-\epsfllx\pspoints
%
%   If `epsfxsize' is 0, we default to the natural size of the picture.
%   Otherwise we scale the graph to be \epsfxsize wide.
%
   \epsfxsize\epsfsize\epsftsize\epsfrsize
   \ifnum\epsfxsize=0 \ifnum\epsfysize=0
      \epsfxsize=\epsftsize \epsfysize=\epsfrsize
%
%   We have a sticky problem here:  TeX doesn't do floating point arithmetic!
%   Our goal is to compute y = rx/t. The following loop does this reasonably
%   fast, with an error of at most about 16 sp (about 1/4000 pt).
% 
     \else\epsftmp=\epsftsize \divide\epsftmp\epsfrsize
       \epsfxsize=\epsfysize \multiply\epsfxsize\epsftmp
       \multiply\epsftmp\epsfrsize \advance\epsftsize-\epsftmp
       \epsftmp=\epsfysize
       \loop \advance\epsftsize\epsftsize \divide\epsftmp 2
       \ifnum\epsftmp>0
          \ifnum\epsftsize<\epsfrsize\else
             \advance\epsftsize-\epsfrsize \advance\epsfxsize\epsftmp \fi
       \repeat
     \fi
   \else\epsftmp=\epsfrsize \divide\epsftmp\epsftsize
     \epsfysize=\epsfxsize \multiply\epsfysize\epsftmp   
     \multiply\epsftmp\epsftsize \advance\epsfrsize-\epsftmp
     \epsftmp=\epsfxsize
     \loop \advance\epsfrsize\epsfrsize \divide\epsftmp 2
     \ifnum\epsftmp>0
        \ifnum\epsfrsize<\epsftsize\else
           \advance\epsfrsize-\epsftsize \advance\epsfysize\epsftmp \fi
     \repeat     
   \fi
%
%  Finally, we make the vbox and stick in a \special that dvips can parse.
%
   \ifepsfverbose\message{#1: width=\the\epsfxsize, height=\the\epsfysize}\fi
   \epsftmp=10\epsfxsize \divide\epsftmp\pspoints
   \vbox to\epsfysize{\vfil\hbox to\epsfxsize{%
      \includegraphics{#1}%
      \hfil}}%
\epsfxsize=0pt\epsfysize=0pt}%
%
%   We still need to define the tricky \epsfaux macro. This requires
%   a couple of magic constants for comparison purposes.
%
{\catcode`\%=12 \global\let\epsfpercent=%\global\def\epsfbblit{%BoundingBox}}%
%
%   So we're ready to check for `%BoundingBox:' and to grab the
%   values if they are found.
%
\long\def\epsfaux#1#2:#3\\{\ifx#1\epsfpercent
   \def\testit{#2}\ifx\testit\epsfbblit
      \epsfgrab #3 . . . \\%
      \epsffileokfalse
      \global\epsfbbfoundtrue
   \fi\else\ifx#1\par\else\epsffileokfalse\fi\fi}%
%
%   Here we grab the values and stuff them in the appropriate definitions.
%
\def\epsfgrab #1 #2 #3 #4 #5\\{%
   \global\def\epsfllx{#1}\ifx\epsfllx\empty
      \epsfgrab #2 #3 #4 #5 .\\\else
   \global\def\epsflly{#2}%
   \global\def\epsfurx{#3}\global\def\epsfury{#4}\fi}%
%
%   We default the epsfsize macro.
%
\def\epsfsize#1#2{\epsfxsize}%
%
%   Finally, another definition for compatibility with older macros.
%

% ================================================================
% execution: why not set
\epsfverbosetrue			% reset at your pleasure
\abovefigskip=\baselineskip		% reset at your pleasure
\belowfigskip=\baselineskip		% reset at your pleasure
\global\let\figtitlefont\bf		% reset at your pleasure
\global\figtitleskip=.5\baselineskip	% reset at your pleasure

\let\nd\noindent %
\font\tenmsb=msbm10
\font\sevenmsb=msbm7
\font\fivemsb=msbm5
\newfam\msbfam
\textfont\msbfam=\tenmsb
\scriptfont\msbfam=\sevenmsb
\scriptscriptfont\msbfam=\fivemsb
\def\Bbb#1{\fam\msbfam\relax#1}

\def\Z{{\Bbb Z}}
\def\R{{\Bbb R}}
\def\C{{\Bbb C}}

\def\a{\alpha}
\def\b{\beta}
\def\l{\lambda}
\def\qed{\hbox{\hskip 6pt\vrule width6pt height7pt depth1pt \hskip1pt}}
\def\vs#1 {\vskip#1truein}
\def\hs#1 {\hskip#1truein}
  \hsize=6truein	\hoffset=.25truein %was \hoffset=1.2truein
  \vsize=8.8truein	%\voffset=1truein
  \pageno=1     \baselineskip=12pt
  \parskip=3 pt		\parindent=20pt
  \overfullrule=0pt	\lineskip=0pt	\lineskiplimit=0pt
  \hbadness=10000 \vbadness=10000 %	REPORT ONLY BEYOND THIS BADNESS
%\nopagenumbers
%\lett

\pageno=1

\footline{\ifnum\pageno=0\hss\else\hss\tenrm\folio\hss\fi}
\hbox{}
\vskip 1truein\centerline{{\bf Some Generalizations of the Pinwheel Tiling}}
\vskip .5truein\centerline{by}
\centerline{Lorenzo Sadun${}^1$}
\footnote{}{1\ Research supported in part by an NSF Mathematical Sciences 
Postdoctoral 
\vs-.1 \hs.15 Fellowship and Texas ARP Grant 003658-037 \hfil}
\vskip .2truein\centerline{Mathematics Department}
\centerline{University of Texas}
\centerline{Austin, TX\ \ 78712}
\vs.1
\centerline{sadun@math.utexas.edu}
\vs.5
\centerline{{\bf Abstract}}
\vs.1 \nd
We introduce a new family of nonperiodic tilings, based on a
substitution rule that generalizes the pinwheel tiling of Conway and
Radin.  In each tiling the tiles are similar to a single triangular
prototile.  In a countable number of cases, the tiles appear in a
finite number of sizes and an infinite number of orientations.  These
tilings generally do not meet full-edge to full-edge, but {\it can} be
forced through local matching rules.  In a countable number of cases,
the tiles appear in a finite number of orientations but an infinite
number of sizes, all within a set range, while in an uncountable
number of cases both the number of sizes and the number of
orientations is infinite.

\vfill\eject 
\nd {\bf \S 1.\ Introduction}
\vs.1
We introduce a new family of nonperiodic tilings, indexed by a
continuous parameter.  In each tiling the tiles are similar to a
single triangular prototile.  In a countable number of cases, the
tiles appear in a finite number of sizes and an infinite number of
orientations.  In a countable number of cases, the tiles appear in a
finite number of orientations but an infinite number of sizes, all
within a set range.  In one case both the number of sizes and
orientations is finite, while in an uncountable number of cases both
the number of sizes and the number of orientations is infinite.  A
piece of a tiling with two sizes and an infinite number of
orientations is shown in Figure 1.

\vs .1
\hs0.95 \vbox{\epsfysize=3truein\epsfbox{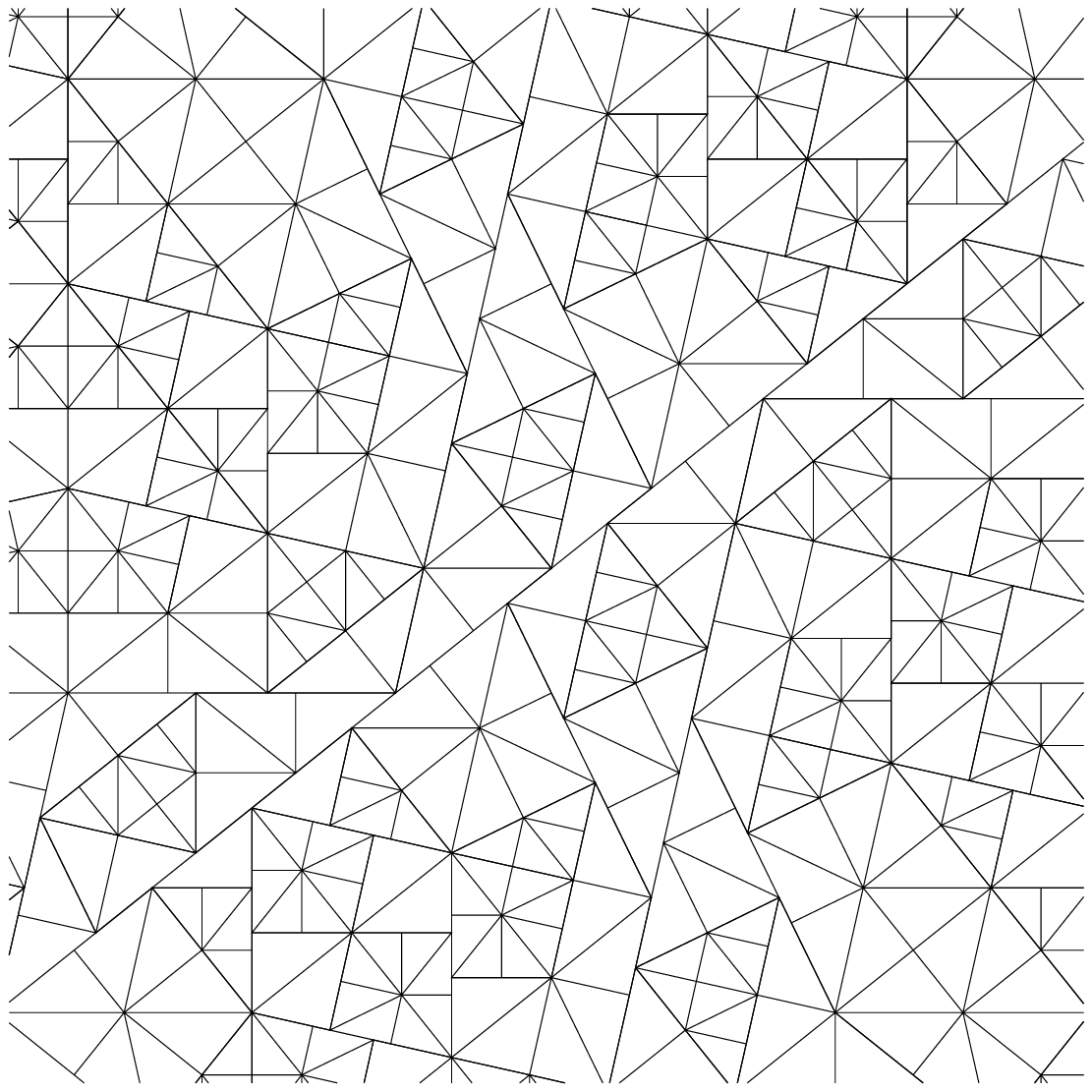}} 
\vs .1
\centerline{Figure 1. Part of the tiling $Til(1/2)$.}
\vs .1

These tilings all arise from a substitution scheme that is quite
similar to the pinwheel tiling of Conway and Radin [R1].  In all cases
the tilings have the ``sibling edge-to-edge" property, which is to say
that two daughter tiles of a single parent tile can only meet
full-edge to full-edge.  However, in general the tiling is not 
{\it globally} edge-to-edge.  Tiles that are not siblings generally do
meet in ways that are not full-edge to full-edge.

Chaim Goodman-Strauss [GS] has recently announced that substitution
tilings with a finite number of prototiles, with each prototile
appearing in a finite number of sizes, and with the sibling
edge-to-edge property, may be associated to local matching rules.
These rules are less restrictive than local atlasses, and in some
cases may allow two adjacent tiles to slide freely along their common edge.
However, the rules are restrictive enough to force any tiling to be
hierarchical (i.e., generated by the substitution scheme), and hence
to be nonperiodic.  In a countable number of cases considered in this
paper, the hypotheses of Goodman-Strauss' theorem are met.

This is striking, because these examples are definitely not globally
edge-to-edge.  Indeed, some of these tilings contain infinitely many
pairs of tiles with edges that partially overlap.  This phenomenon,
called ``slippage'', can be seen in Figure 1, especially along the
long diagonal.  The lengths of the overlaps are not the same for all
pairs of tiles.  In the tiling of Figure 1, in fact, overlaps occur
with an {\it infinite} number of distinct lengths.  Thus we have the
(seemingly) paradoxical situation that Goodman-Strauss' local matching
rules force a hierarchical pattern that, in turn, forces an infinite
variety of local behavior.

The tilings of this paper also suggest a relaxation of the rules of
the tiling dynamical systems game.  Prior to Conway and Radin's work,
a tiling dynamical system required that all tiles be generated from a
finite set of prototiles by translation only, with rotation and
reflection not allowed.  By those rules, the pinwheel tiling, which
involves an infinite number of orientations of a single triangle,
uniformly distributed about the circle, was not a tiling!  The pinwheel
example helped force a reconsideration of the rules, and a new
consensus has emerged that tiles may be generated from prototiles by
Euclidean motions, not just by translations.

The present examples suggest a further extension from the Euclidean
group to the conformal group.  In almost all cases, these tilings
consist of a single prototile appearing in an infinite number of sizes
and an infinite number of orientations.  Although the distribution of
sizes is not constant (which it could not be, as the conformal group
has no Haar measure), it is described by a piecewise-constant
function.

Such symmetry is related to properties of the spectrum.  In the case
of the pinwheel, Radin [R2] showed that the statistical rotational
invariance of the system results in a spectrum that is rotationally
invariant, and hence contains no discrete component.  However, Radin's
argument did not address the dependence of the spectrum on the radial
variable, which corresponds to a length scale.  Almost all the
examples in this paper have a joint distribution of sizes and
orientation that is rotationally invariant, proving (by Radin's
argument) that they have continuous spectrum.  In addition, in almost
all these examples the distribution is absolutely continuous in the
size parameter.  We conjecture that, in these cases, the spectrum is
purely absolutely continuous.

Purists who object to the extension to the conformal group may
restrict their attention to the countable and dense set of examples
for which only a finite number of sizes appear.  They may then take
these sizes of the basic triangle as their prototiles, which then
generate all tiles via Euclidean motions.  We call such examples
``rational tilings''.

The division of the paper is as follows.  In \S 2 we explain the
substitution rules and the resulting constructions of the tilings, and
classify right triangles according to whether they generate finite
numbers of sizes and orientations, finite sizes and infinite
orientations, infinite sizes and finite orientations, or infinite
sizes and infinite orientations.  In \S 3 we examine in detail an
example in which tiles appear in two sizes and an infinite number of
orientations.  We prove that this tiling is sibling edge-to-edge but
exhibits slippage; this tiling contains pairs of tiles that meet in an
infinite number of distinct ways.  In \S 4 we compute the limiting
distribution of sizes and orientations in all rational tilings.  In \S
5 we compute the statistical distribution of sizes and orientations in
irrational tilings.  In \S 6 we examine two exceptional rational
tilings.

\vs.2
%\vfill\eject
\nd {\bf \S 2.\ Definitions and constructions}
\vs.1

As with any substitution scheme, we obtain a tiling of the plane by a
succession of subdivisions (``deflations''), rescalings and
repositionings.  We start with a single prototile $T$, define a
subdivision rule, and let $T_n$ be the result of subdividing the basic
tile $n$ times.  We will construct infinite sequences of integers $N_1
< N_2 < \cdots$, rescalings $s_i$ and Euclidean motions $e_i$ such
that $e_i(s_i(T_{N_i}))$ is a proper subtiling of
$e_{i+1}(s_{i+1}(T_{N_{i+1}}))$.  By taking the union of the finite
tilings $e_i(s_i(T_{N_i}))$, we obtain a tiling of an infinite region,
typically the entire plane.

Our task is complicated by the fact that there is not a
straightforward connection between $N_i$ and $s_i$.  (Unlike, say, the
pinwheel tiling, where $s_i$ is an expansion by a linear factor
$5^{N_i/2}$.) We cannot, in general, define a simple rule for
deflating and rescaling by a fixed factor, and simply iterate this
rule.  We must be more subtle.  However, in the end, everything does
work out.
 
Our basic prototile is a right triangle, with base $b$, altitude $a$
and hypotenuse $c$, as shown in Figure 2.  We refer to the base as the
``long'' leg, although in principle we might have $b \le a$. At this
point we put no restrictions on the acute angle $\theta =
\sin^{-1}(a/c)$.  We divide $T$ into five triangles as in Figure 3.
Subtriangles 1--4 are similar to $T$, but are smaller by a linear
factor of $b/2c$.  Subtriangle 5 is also similar to $T$, but is
smaller than $T$ by a linear factor of $a/c$.  This completes the
first subdivision of $T$.  In the $n$-th subdivision of $T$, we take
the largest triangles in $T_{n-1}$, and subdivide each of them into
five similar triangles.

\vs .1
\hs1.15 \quad \vbox{\epsfxsize=3truein\epsfbox{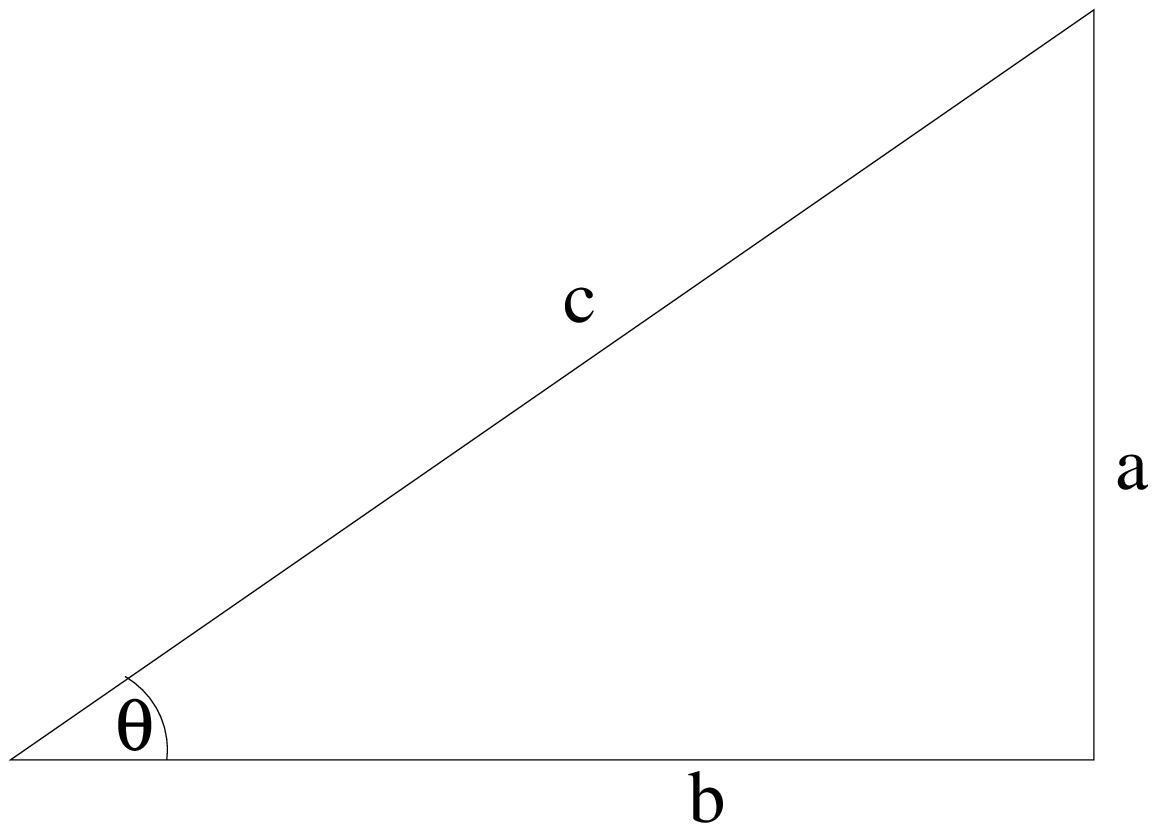}} 
\vs .1
\centerline{Figure 2. The basic triangle.}
\vs .1

The pattern of subdivision depends on the angle $\theta$.  For
example, if $a > b/2$, then the largest triangle in $T_1$ has
hypotenuse $a$, and only this triangle is subdivided in step 2.  Thus
in $T_2$ there are 9 triangles, one with hypotenuse $a^2/c$, 4 with
hypotenuse $ab/2c$, and 4 with hypotenuse $b/2$, as shown in Figure 4.
The third subdivision depends on whether $b/2 > a^2/c$, in which case
the 4 triangles of hypotenuse $b/2$ are subdivided next, or whether
$a^2/c>b/2$, in which case the one triangle of hypotenuse $a^2/c$ is
subdivided once again.

\vs .1
\hs1.15 \vbox{\epsfxsize=3truein\epsfbox{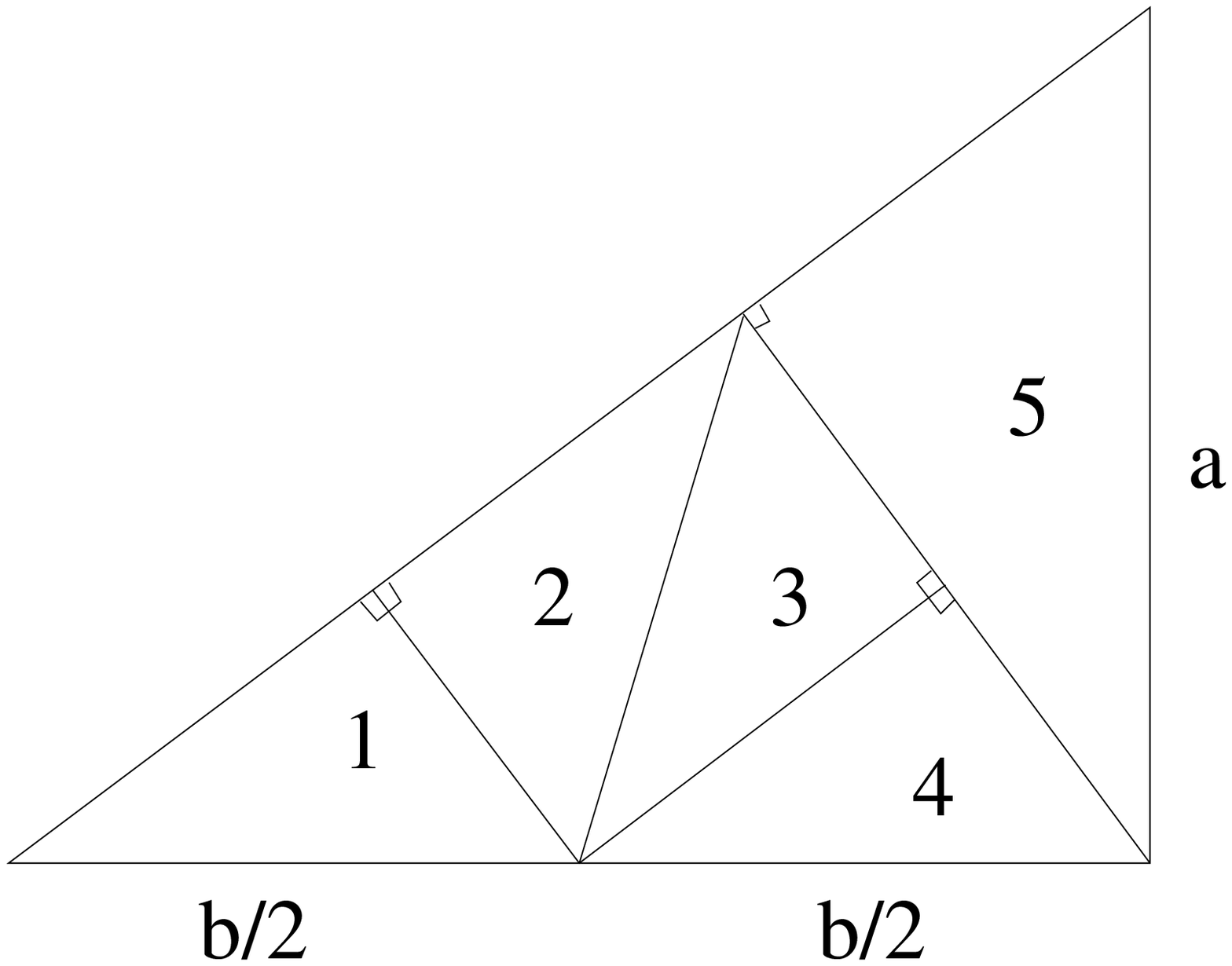}} 
\vs .1
\centerline{Figure 3. The substitution rule.}
\vs .1

\vs .1
\hs1.15 \vbox{\epsfxsize=3truein\epsfbox{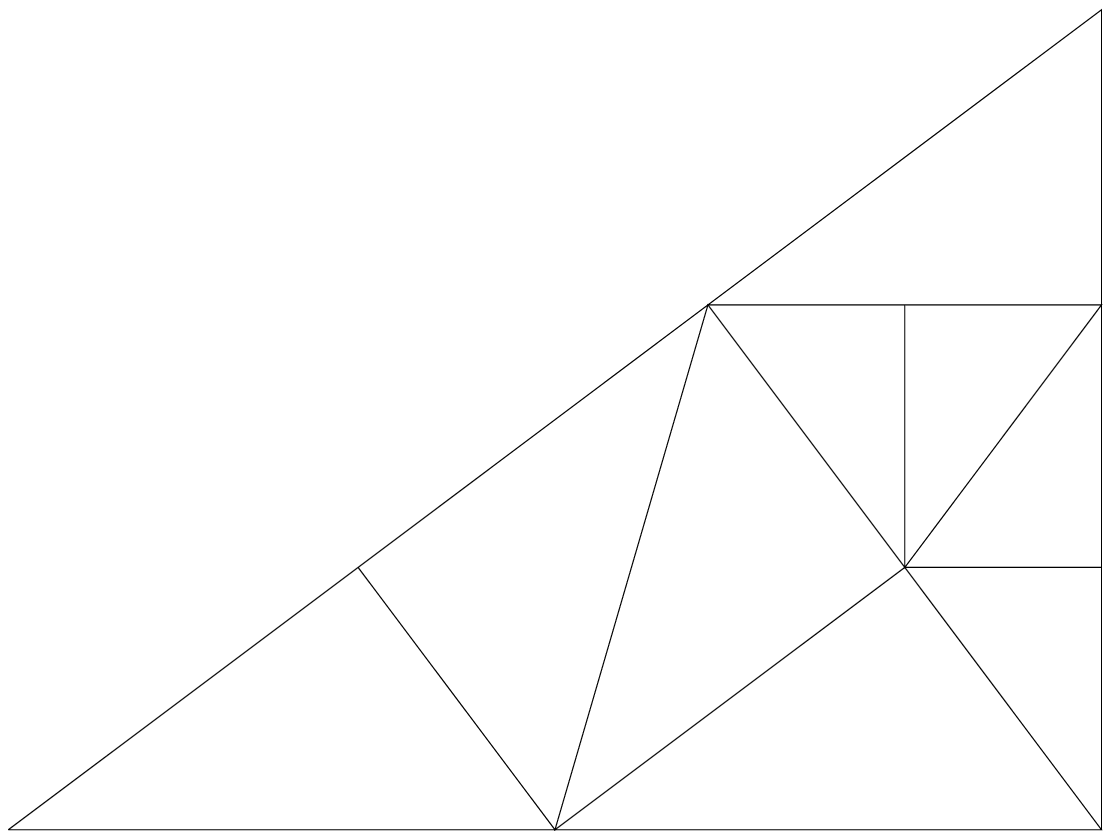}} 
\vs .1
\centerline{Figure 4. Two subdivisions of a triangle with $a>b/2$.}
\vs .1

\vs .1
\hs1.15 \vbox{\epsfxsize=3truein\epsfbox{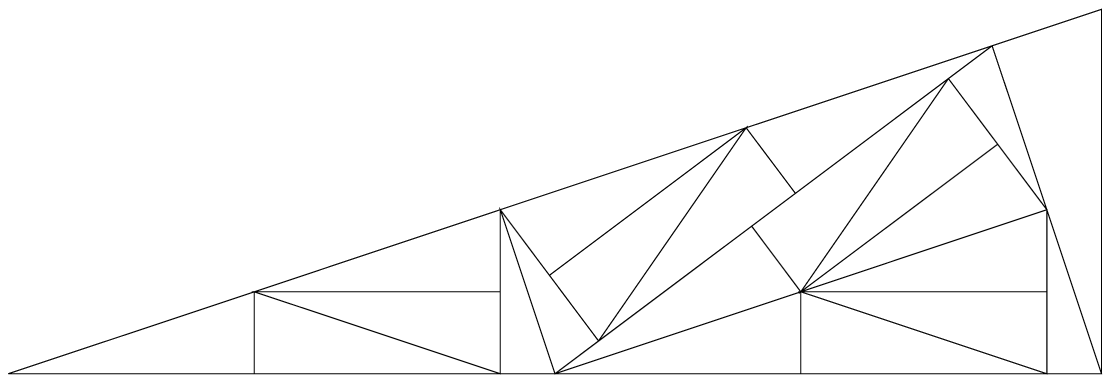}} 
\vs .1
\centerline{Figure 5. Two subdivisions of a triangle with $a<b/2$.}
\vs .1

The subdivision of a triangle with $b/2 > a$ is shown in Figure 5.
$T_2$ has 21 triangles, 16 of hypotenuse $b^2/4c$, 4 of hypotenuse
$ab/2c$, and one of hypotenuse $a$.

Finally, if $b=2a$, then all five triangles in $T_1$ have the same
size.  All five are subdivided in the next stage, yielding $T_2$ with
25 congruent triangles, all of which are then subdivided to give $T_3$
with 125 congruent triangles, and so on.  This is the pinwheel tiling
of Conway and Radin [R1].

We will show that, for any angle $\theta$, this subdivision scheme
generates nonperiodic tilings of the plane.  We first need three
technical lemmas:

{\nd \bf Lemma 1:} {\it In $T_n$, the ratio of the hypotenuse of the largest
tile to the hypotenuse of the smallest tile is strictly less than 
$\max(c/a,2c/b)$.}

\nd {\it Proof:} The proof is by induction on $n$.  It is clearly true
for $T_1$.  Now suppose that in $T_{k}$ the largest tile has
hypotenuse $H$, while the smallest tile has hypotenuse $h$, with $H/h
< \max(c/a,2c/b)$.  After subdividing all the triangles with
hypotenuse $H$, the largest hypotenuse is strictly less than $H$,
while the smallest hypotenuse is the smallest of $h$, $aH/c$, and
$bH/2c$. The ratio of largest to smallest hypotenuses in $T_{k+1}$ is
therefore strictly less than the largest of $H/h$, $c/a$ and
$2c/b$. \qed

We will often consider an individual tile in $T_n$, and examine the
effect of further subdivision of $T_n$ on that tile.  When we write $t
\in T_n$, we mean that $t$ is a single tile.  When we write $t \subset
T_N$, for some $N>n$, we mean the collection of tiles in $T_N$ whose
union is $t$.  When we write that $t \subset T_N$ is similar to $T_k$,
we mean that the individual tiles in this collection fit together to
form the single tile $t \in T_n$ in exactly the same pattern that the
individual tiles of $T_k$ fit together to form $T$.

{\nd \bf Lemma 2:} {\it Given an integer $n$ and a tile $t \in T_n$,
there exists an integer $N>n$ such that $t \subset T_N$ is the union
of more than one tile.  That is, every tile in $T_n$, no matter how
small, is eventually subdivided further.}

\nd {\it Proof:} The number of tiles grows without bound (increasing by at
least 4 every turn).  By Lemma 1, the area of each tile is at least
$\min(a^2/c^2, b^2/4c^2)$ times the area of the largest tile.  Thus,
as $N\to \infty$, the area of the largest tile in $T_N$ must go to
zero.  For $N$ large enough, this maximal area will be less than the
area of $t$, indicating that the tile $t$ has necessarily been
subdivided. \qed

{\nd \bf Lemma 3:} 
{\it Given an integer $n$, a tile $t \in T_n$ and an integer $m$, 
there exists an integer $N \ge n$ such that 
$t \subset T_N$ is similar to $T_m$.}

\nd {\it Proof:} The proof is by induction on $m$.  First we consider
$m=1$.  Let $N$ be the smallest integer such that $t \in T_N$ is
subdivided.  By Lemma 2, such an $N$ exists.  Since $t \subset
T_{N-1}$ is a single tile, $t \subset T_N$ is the result of taking a
single triangle and applying the deflation rule of Figure 3 once.
Thus $t \in T_N$ is similar to $T_1$.

Now suppose the theorem has been proved for $m=k$, and suppose that $t
\subset T_{N_0}$ is similar to $T_k$.  We consider what happens at the
$N_0+1$st step.  The largest triangles in $T$ get subdivided.  This
means that either the largest triangles in $t$ get subdivided, or that
none of the triangles in $t$ get subdivided.  In the first case $t$
will then become similar to $T_{k+1}$, and we are done.  In the second
case we consider the possibilities for the $N_0+2$nd step, and so on.
By Lemma 2, the largest triangle in $t$ must {\it eventually} be
subdivided, say on the $N$-th turn, so $t \subset T_N$ will be similar
to $T_{k+1}$.  \qed

{\nd \bf Definition:} A {\it supertile of order $n$} is a collection of
tiles that is similar to $T_n$.  Equivalently, a supertile of order $n$
is a region of the form $e(s(T_n))$, where $s$ is a rescaling and $e$ is 
a Euclidean motion.

%\vfill\eject
 
{\nd \bf Theorem 1:} {\it Given any right triangle $T$, there exist
tilings of the plane with right triangles similar to $T$, such that
any finite set of tiles lies in a supertile, and such that the areas
of tiles are bounded both above and below.}
 
{\nd \it Proof:} First pick a succession of integers $n_1, n_2,
\ldots$ and tiles $t_i \in T_{n_i}$.  Let $N_1=n_1$.  Pick additional
integers $N_i$, $i=2,3,\ldots$, such that, taking $t_i \in T_{n_i}$
and subdividing $T_{n_i}$ an additional $N_i-n_i$ times, $t_i \subset
T_{N_i}$ is similar to $T_{N_{i-1}}$.  By Lemma 3, such integers
always exist.

Now pick a triangle similar to $T_0$ and place it in the plane.  One can
construct a supertile $S_1$ of order $N_1$ that contains this triangle in
position $t_1$.  One then constructs a supertile $S_2$ of order $N_2$ 
such that $S_1$ sits inside $S_2$ as $t_2$ sits inside $T_{n_2}$.  One
continues the process, building supertile $S_{k+1}$ such that $S_k$ sits
inside $S_{k+1}$ as $t_{k+1}$ sits inside $T_{n_{k+1}}$.

The union of all the supertiles is a tiling of an infinite region.
For almost all choices of the $t_i$'s (e.g., having the edges of $t_i$
lie in the interior of $T_{n_i}$ infinitely often), this region will
be the entire plane.

Since the ratio of largest to smallest triangle is uniformly bounded
for $S_n$ by Lemma 1, no tile may have hypotenuse longer than
$\max(c/a, 2c/b)$ times the hypotenuse of a fixed tile, and no tile
may have hypotenuse less than $\min(a/c,b/2c)$ times the hypotenuse of
the same fixed tile.  This provides both an upper and lower bound to
the size of the tiles. \qed

In this construction, many choices were made.  Different choices generally
lead to different tilings, but these different tilings have many
properties in common.  The following Theorems 2 and 3 apply to all 
tilings constructed in the manner of Theorem 1.

{\nd \bf Theorem 2:} {\it In a tiling, the number of orientations in
which the basic triangle appears is finite if $\theta/\pi$ is rational
and infinite if $\theta/\pi$ is irrational.}

\nd {\it Proof:} First suppose that $\theta/\pi$ is rational. Consider
a tile, positioned as in Figure 3, with side $b$ along the $x$ axis.
Let $P$ denote reflection about the $x$ axis, and let $R_\alpha$
denote a counterclockwise rotation by angle $\alpha$.  After
subdividing once, the orientations of the five daughter tiles,
relative to the parent tile, are given by the following elements of
$O(2)$: $R_\theta P$ (twice), $R_\theta$, $R_{\pi + \theta}$, and
$R_{\pi/2 + \theta} P$.  The orientations of tiles in a further
subdivision are words in these five elements of $O(2)$.  However, with
$\theta/\pi$ rational, these five elements generate a finite subgroup
of $O(2)$, so only a finite number of orientations can ever appear in
a future subdivision.

Since any region of our tiling of the plane sits inside a supertile,
any two tiles must have their orientations, relative to the supertile
itself, in this group.  Thus their orientations, relative to each
other, and hence to a fixed reference tile, must lie in the group.
Thus only a finite number of orientations can appear in the tiling.

Now suppose that $\theta/\pi$ is irrational.  In a basic subdivision,
we will keep track only of the 4 triangles of hypotenuse $b/2$,
ignoring the triangle of hypotenuse $a$.  When these 4 triangles
divide, we will only keep track of the 16 resultant triangles of side
$b^2/4c$, and so on.  In the second generation we find orientations $1
= (R_\theta P)^2$ and $R_{2 \theta}=R_\theta^2$, among others.  In the
$2n$-th generation we find $1, R_{2 \theta}, R_{4 \theta}, \ldots,
R_{2 n \theta}$.  Since $\theta/\pi$ is irrational, these $2n+1$
orientations are distinct.  Since our tiling of the plane contains
supertiles of arbitrarily large size, there is no bound to the number
of different orientations that appear. \qed
    
{\nd \bf Theorem 3:} {\it In a tiling, the number of sizes in which
the basic triangle appears is infinite if
$\ln(\sin(\theta))/\ln[\cos(\theta)/2]$ is irrational and finite if
$\ln(\sin(\theta))/\ln[\cos(\theta)/2]$ is rational.  In particular,
if $\ln(\sin(\theta))/\ln[\cos(\theta)/2]=p/q$, with $p$ and $q$
relatively prime integers, then the number of sizes in the tiling is
$\max(p,q)$.}

\nd {\it Proof:} Let $A = a/c$ and let $B=b/2c$.  If $\ln(A)/\ln(B)$ is
irrational, the only way two monomials $A^aB^b$ and $A^cB^d$ can equal
is if $a=c$ and $b=d$.  We will show that the sizes of triangles in
$T_n$ (relative to the original triangle) is given by such monomials,
and that the number of distinct powers of $A$ grows without bounds as
$n \to \infty$. This will show that the number of distinct sizes grows
without bound a $n \to \infty$.

In each subdivision there are 4 tiles of size $B$ relative to the
parent and one tile of size $A$.  Thus the descendants of a given tile
all have sizes that are monomials $A^aB^b$ relative to the ancestor.
For every $n>0$, $T_n$ contains at least one tile with size $A^0B^b$;
just take a $B$ child of a $B$ child of $\ldots$ of one of the
original $B$ children (or the $B$ child itself, if it has not
subdivided).  For every $n>0$, $T_n$ contains at least one tile with
size $A^1B^b$; take a $BBB\ldots$ descendant of the original $A$
child.  Once $n$ is large enough to have the original $A$ child
divide, there is at least one tile with size $A^2 B^b$.  In general,
once $n$ is large enough to allow a $a$-th generation $AA\ldots$
child, it will always have at least one tile whose size has exactly
$a$ powers of $A$.  This completes the irrational case.

Now suppose that $\ln(A)/\ln(B)=p/q$, with $p$ and $q$ relatively
prime.  Thus $A^q = B^p$.  Let $r = A^{1/p} = B^{1/q}$.  Every
monomial $A^aB^b$ is a power of $r$.  Assume for the moment that $p
\le q$.  By Lemma 1, the ratio of sizes of any two tiles is greater
than $B=r^q$.  Thus only at most $q$ distinct sizes can appear.  To
see that $q$ sizes {\it do} appear, we note that $A^0B^{b_0}, A^1
B^{b_1},
\ldots, A^{q-1} B^{b_{q-1}}$ are all distinct powers of $r$.  

If $p > q$, Lemma 1 states that the ratio of any two sizes is at least
$A=r^p$, so at most $p$ different sizes can occur.  We produce $p$
different sizes by examining different powers of $B$.  In either case,
the number of distinct sizes is $\max(p,q)$. \qed

We refer to tilings with $\ln(A)/\ln(B)=p/q$ as {\it $(p/q)$ rational
tilings}, and denote the class of such tilings as $Til(p/q)$.  For $z$
irrational, we will similarly denote the class of tilings with
$\ln(A)/\ln(B)=z$ as $Til(z)$.  The different tilings in a class are
all derived from the same substitution rule, and have many properties
derivable from this rule.  When discussing such properties, we will
sometimes refer to a typical element of the class as ``the tiling
$Til(z)$''.

Note that $\ln(\sin(\theta))/\ln[\cos(\theta)/2]$ is a strictly
decreasing function of $\theta$ on the interval $(0,\pi/2)$.  From
this monotonicity, and from the countability of the rationals, it is
clear that only a countable set of angles $\theta$ give rise to a
finite number of rotations, and only a countable set of angles
$\theta$ give rise to a finite number of sizes.  The intersection of
these two countable sets turns out to be a single point.

{\nd \bf Theorem 4:} {\it  The only angle that gives rise to 
both a finite number of orientations and a finite number of 
sizes is $\theta=\pi/4$.  That is, the tiling $Til(1/3)$.}

\nd {\it Proof:} Let $x = \exp(i \theta)$, with $0<\theta<\pi/2$.  
We are looking for solutions to the equation
$\sin(\theta)^q = [\cos(\theta)/2]^p$, which we rewrite as
$$ 
2^q (x + {\bar x})^p = 2^{2p}(-i)^{q}(x-{\bar x})^q, \eqno(2.1) 
$$
where ${\bar x}=\exp(-i\theta)=x^{-1}$. Note that, for
fixed $p$, $q$, there
is at most one solution to equation (1) in the first quadrant, since
$\ln[\sin(\theta)]/\ln[\cos(\theta)/2]$ is monotonic.

Since by assumption $\theta$ is a rational multiple of $\pi$, $x$ is a
primitive $n$-th root of unity for some integer $n$.  If $x$ is a
solution and $q$ is even, then equation (1) has (real) integer
coefficients, and {\it all} the primitive $n$-th roots of unity are
also solutions.  If $q$ is odd, all the primitive $n$-th roots of
unity are solutions either to equation (1) or to the conjugate
equation
$$ 
2^q (x + {\bar x})^p = 2^{2p}i^{q}(x-{\bar x})^q. \eqno(2.2) 
$$
Equation (2), with $q$ odd, has no solutions in the first quadrant, as
the right hand side is positive but the left hand side is negative.
Since equation (1) admits only one solution in the first quadrant,
there must be exactly one primitive $n$-th root of unity in the first
quadrant.

This means that $n$ must equal 5, 6, 7, 8, 10, 12 or 18.  Checking
these individually, we see that only $n=8$, or $\theta=\pi/4$, yields
a rational value of $\ln[\sin(\theta)]/\ln[\cos(\theta)/2]$.  In that
one case $\sin(\theta)=\cos(\theta)=\sqrt{2}/2$, and $\sin^3(\theta) =
\cos(\theta)/2$.  \qed

In a periodic tiling, all of the sizes and orientations are exhibited
in a compact region, so neither the number of sizes nor the number
of orientations can be infinite.  Thus we have

{\bf \nd Corollary} {\it If $z \ne 1/3$, then the tiling $Til(z)$ is
not periodic}.  

In fact, it will turn out that $Til(1/3)$ is not periodic, either.
This will be shown in \S 6.

{\nd \bf \S  3.  An example with two sizes.}

In this section we consider in detail the tiling $Til(1/2)$, shown in
Figure 1.  This example is chosen not as a special case, but rather as
a simple example of some general phenomena.  Based on the statistical
analysis of \S 4, we expect all tilings $Til(p/q)$, with $q>1$ and
$p/q \ne 1/3$, to be qualitatively similar to $Til(1/2)$.
Specifically, in all these cases the population matrix has two or more
eigenvalues with modulus bigger than one.  This causes fluctuations in
the statistical composition of supertiles $S_n$ to grow with $n$.
This, in turn, can cause the tiling to fail to be globally
edge-to-edge.

\vfill\eject

{\nd \bf Theorem 5:} {\it $Til(1/2)$ is a tiling with two sizes of
tiles, each of which appears in an infinite number of orientations.
The substitution scheme has the sibling edge-to-edge property, but the
tiling is not globally edge-to-edge.  Specifically, the tiles meet in
an infinite number of ways.}

\nd {\it Proof:} The existence of the tiling, the number of sizes and
the number of orientations follow from Theorems 1--3.  The sibling
edge-to-edge property is manifest, if we consider the basic tile to
have four vertices -- the three obvious ones and the midpoint of the
long leg.  The difficulty is in proving that tiles meet in an infinite
number of distinct ways.

\vs .1
\hs0.65 \vbox{\epsfxsize=4truein\epsfbox{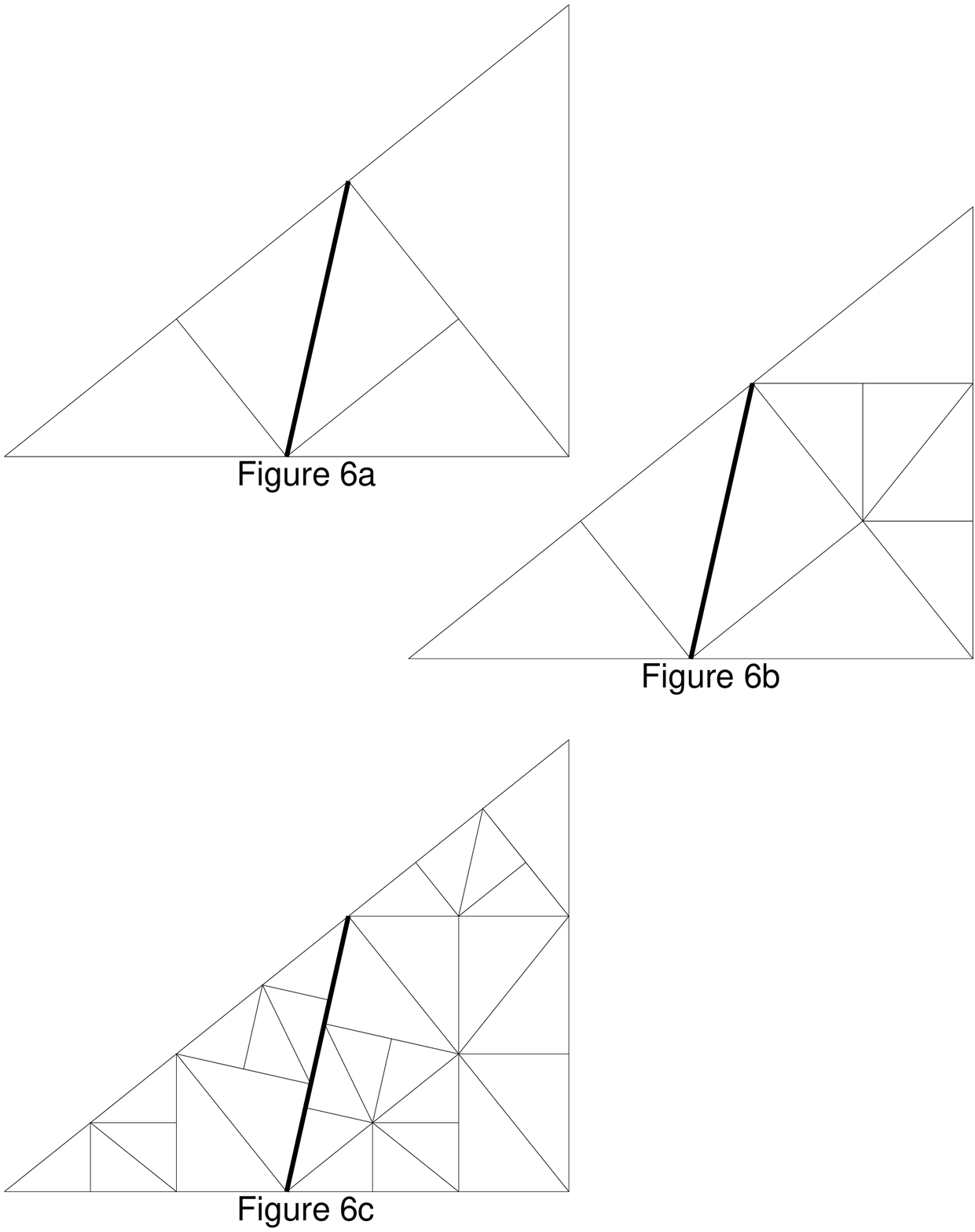}} 
\vs .2
\centerline{Figure 6. Three stages of subdivision for $Til(1/2)$.}
\vs .1

The process by which this happens is as follows.  There are certain
special lines in the tiling.  The long diagonal in Figure 1 is an
example.  On each side of such a line there are triangles, all of the
same size, whose hypotenuses or long legs make up part of the long
line.  The pattern is {\it different} on the two sides of the line,
with one side having (say) more legs in a certain region and the other
side having more hypotenuses.  This imbalance causes the tiles on one
side of the line to appear shifted relative to those on the other
side, a phenomenon we call ``slippage''.  We will exhibit regions where the
imbalance is arbitrarily large.  Since the imbalance is unbounded, the
slippage of one side relative to the other reaches arbitrarily high
multiples of $b$, modulo $c$, where $b$ and $c$ are the lengths of the
long leg and hypotenuse.  Since $b/c$ is irrational, this gives rise
to an infinite number of ways in which one triangle can meet another
across the long line.

Lines where slippage occurs are called ``fault lines''.  Note that
each fault line has only finite length, and allows only a finite
amount of slippage.  However, we will find fault lines of arbitrarily
long length with arbitrarily much slippage.  This precludes there
being a finite bound on how many ways one triangle can meet another.

The $Til(1/2)$ tiling is based on a right triangle with legs
$a=\sqrt{2(\sqrt{17}-1)}$ and $b=\sqrt{17}-1$ and hypotenuse $c=4$.
Several degrees of subdivision are shown in Figure 6.  At each level
there are two sizes of triangle, whose linear sizes differ by a factor
of $a/c$.  The heavily shaded lines in Figure 6 are fault lines.  For
each $T_n$, with $n$ even, only large triangles abut the illustrated
fault line.
%Smaller fault lines appear after additional subdivisions.
%In particular, by looking at the descendants of the big triangle of $T_1$,
%one can find a scale copy of $T_{n-1}$ inside $T_n$.   
Once a fault line is formed, the triangles on opposite
sides of the line evolve separately, and begin to slip.  We shall
prove that this slippage increases without bound.

$Til(1/2)$ may be viewed as a traditional substitution system, with
two prototiles, which we call $B$ and $S$ (for big and small).  Each
subdivision, followed by rescaling by $c/a$, may be viewed as a
replacement of each $S$ triangle by a $B$ triangle, and replacement of
each $B$ triangle by a $B$ triangle and four $S$ triangles.  Let the
population of $T_n$ be $\Psi_n = \left (\matrix{N_B \cr N_S}
\right )$, where $N_B$ and $N_S$ are the numbers of big and small tiles.
$\Psi$ satisfies $\Psi_{n+1} = M \Psi_n$, where the population matrix is
$$ M= \left (\matrix{1 & 1 \cr 4 & 0} \right ), \eqno (3.1) $$
with eigenvalues $\lambda_{\pm} = (1 \pm \sqrt{17})/2$ and eigenvectors
$\zeta_{\pm} = \left ( \matrix {\sqrt{17} \pm 1 \cr \pm 8} \right )$.
As $n$ grows, the ratio of $N_B$ to $N_S$ approaches $(\sqrt{17} + 1)/8
\approx 0.6404$.  The exact populations are
$$ \eqalign{ N_B(n) = & 
{1 \over \sqrt{17}} \left ( \left ({\sqrt{17} + 1 \over 2}
\right )^{n+1} -  \left ({1 - \sqrt{17} \over 2}
\right )^{n+1}\right ) \cr
N_S(n) = 4 N_L(n-1)= &
{4 \over \sqrt{17}} \left ( \left ({\sqrt{17} + 1 \over 2}
\right )^{n} -  \left ({1 - \sqrt{17} \over 2}
\right )^{n}\right ).}  \eqno(3.2) $$
Note that $|\lambda_-|>1$, so that $|8 N_B - (\sqrt{17}+1)N_S|$ grows
with $n$.

Next we consider what happens along a fault line.  To do this we must
consider the boundary of $T_{2n}$.  Note that the hypotenuse and long
leg of $T_2$ consist only of hypotenuses and long legs of big
triangles.  Applying the subdivision again, we get that the hypotenuse
and long leg of $T_4$ also consists only of hypotenuses and long legs
of big triangles.  Similarly for all $T_{2n}$.

The evolution of these legs and hypotenuses is a substitution system
in its own right, only in one dimension.  There are four symbols,
$H^+$, $H^-$, $L^+$ and $L^-$, representing the two orientations of
the hypotenuse and long leg, respectively.  From Figure 6b, we see
that the substitution rule, which we denote $\sigma_0$, is
$$ 
\sigma_0(H^+) = L^+ L^- H^+, \qquad \sigma_0(H^-) = H^- L^+ L^- ,
\qquad \sigma_0(L^\pm)= H^\pm H^\pm. \eqno(3.3) 
$$
Since $L^+$ and $L^-$ only appear in the combination $L^+ L^-$, we
can define a new symbol $L=L^+ L^-$ and have a substitution system
with 3 elements, whose rule we denote $\sigma$:
$$
\sigma(H^+)= L H^+, \qquad \sigma(H^-)= H^- L, \qquad 
\sigma(L)= H^+ H^+H^- H^-.
\eqno (3.4) $$

{\nd \bf Lemma 4:} {\it The sequence $\sigma^n(H^+)$ 
contains neither the subsequence $LL$ nor the subsequence $H^-H^+$.}

\nd {\it Proof:} The proof is simple induction.  The only way to generate 
an $LL$ is from an $H^-H^+$, and the only way to generate an $H^-H^+$
is from $LL$.  Since neither appear in the first generation, neither
appears in any subsequent generation. \qed

{\nd \bf Lemma 5:} {\it The sequence $\sigma^n(H^+)$ does not contain
a subsequence of more than 6 consecutive $H$'s.}

\nd {\it Proof:} Since $LL$ does not occur in $\sigma^{n-1}(H^+)$, 
the longest possible sequence of $H$'s in $\sigma^n(H^+)$ would come
from a sequence $H^+ L H^-$ in $\sigma^{n-1}(H^+)$.  This gives rise
to $L H^+H^+H^+H^-H^-H^-L$, or 6 $H$'s in a row.  \qed

Let $f(n)$ equal the number of $L$'s in the first
half of the sequence $\sigma^n(H^+)$ minus the number of $L$'s in the 
second half of the sequence.  As we shall see, $f(n)$ is closely related
to the extent to which slippage occurs along the largest fault line
in $T_{2n+2}$.

{\nd \bf Lemma 6:} {\it If $|f(n)|>6$, then $|f(n+1)| \ge |f(n)| +2$.}

Let $s$ and $s'$ denote the first and second halves of $\sigma^n(H^+)$,
respectively.  Suppose that $s$ contains $h$ $H$'s and
$l$ $L$'s, while $s'$ contains $h'$ $H$'s and
$l'$ $L$'s.  Note that $f(n)=l-l'=h'-h$.  Since each $H$ generates an $H$ and
an $L$, while each $L$ generates four $H$'s, $\sigma(s)$ contains 
$4l+h$ $H$'s and $h$ $L$'s, while $\sigma(s')$ contains 
$4l'+h'$ $H$'s and $h'$ $L$'s.  Thus $\sigma(s)$ contains $2f(n)$ more
terms, but $f(n)$ fewer $L$'s, than $\sigma(s')$.

Now suppose $f(n) > 0$. The first half of $\sigma^{n+1}(H^+)$ is all
of $\sigma(s)$, minus the last $f(n)$ elements, while the second half
of $\sigma^{n+1}(H^+)$ is the last $f(n)$ elements of $\sigma(s)$ and
all of $\sigma(s')$.  Thus $f(n+1)$ equals $-f(n)$ minus twice the
number of $L$'s in the last $f(n)$ elements of $\sigma(s)$.  Since
$f(n)>6$, there must be at least one $L$ in the last $f(n)$ elements
of $\sigma(s)$, so $|f(n+1)| = -f(n+1) \ge f(n)+2$.

If $f(n)<0$, then the first half of $\sigma^{n+1}(H^+)$ is all of
$\sigma(s)$, plus the first $|f(n)|$ elements of $\sigma(s')$.  We
then have $f(n+1)$ equalling $-f(n)$ plus twice the number of $L$'s in
the first $|f(n)|$ elements of $\sigma(s)$.  By Lemma 5, there must be
at least one such $L$, so $|f(n+1)=f(n+1) \ge 2-f(n)= |f(n)|+2$.
\qed

{\nd \bf Lemma 7:} {\it $\lim_{n \to \infty} |f(n)| = +\infty$}.

\nd {\it Proof:} By explicit computation, $f(1)=1$, $f(2)=-1$,
$f(3)=1$, $f(4)=-3$, $f(5)=3$, $f(6)=-5$ and $f(7)=9$.  By Lemma 6,
for $n \ge 7$, $|f(n+1)| >  |f(n)|$, so $|f(n)| \ge n+2$ goes to infinity
as $n \to \infty$. \qed

We have proven that $|f(n)|$ grows without bound, which is all that we
need.  In fact, $|f(n)|$ grows exponentially.  For large $n$, an
approximate fraction $2/(3+\sqrt{17})$ of the elements of
$\sigma^n(H^+)$ are $L$'s, so $f(n+1) \approx -f(n) + 2 f(n)
\times 2/(3+\sqrt{17}) = f(n) (1 - \sqrt{17})/2$.  The growth rate,
$(1 - \sqrt{17})/2$, equals $\lambda_-$, the second eigenvalue of
the population matrix $M$.

We now return to the question of slippage along fault lines.  Consider
two large triangles that meet hypotenuse to hypotenuse to form a
rectangle, as in Figure 7.  Let $P$ and $R$ be the ends of the common
hypotenuse, and let $Q$ be the midpoint.  Rotation by $\pi$ about $Q$
sends each triangle into the other.

\vs .1
\hs1.70 \vbox{\epsfxsize=2truein\epsfbox{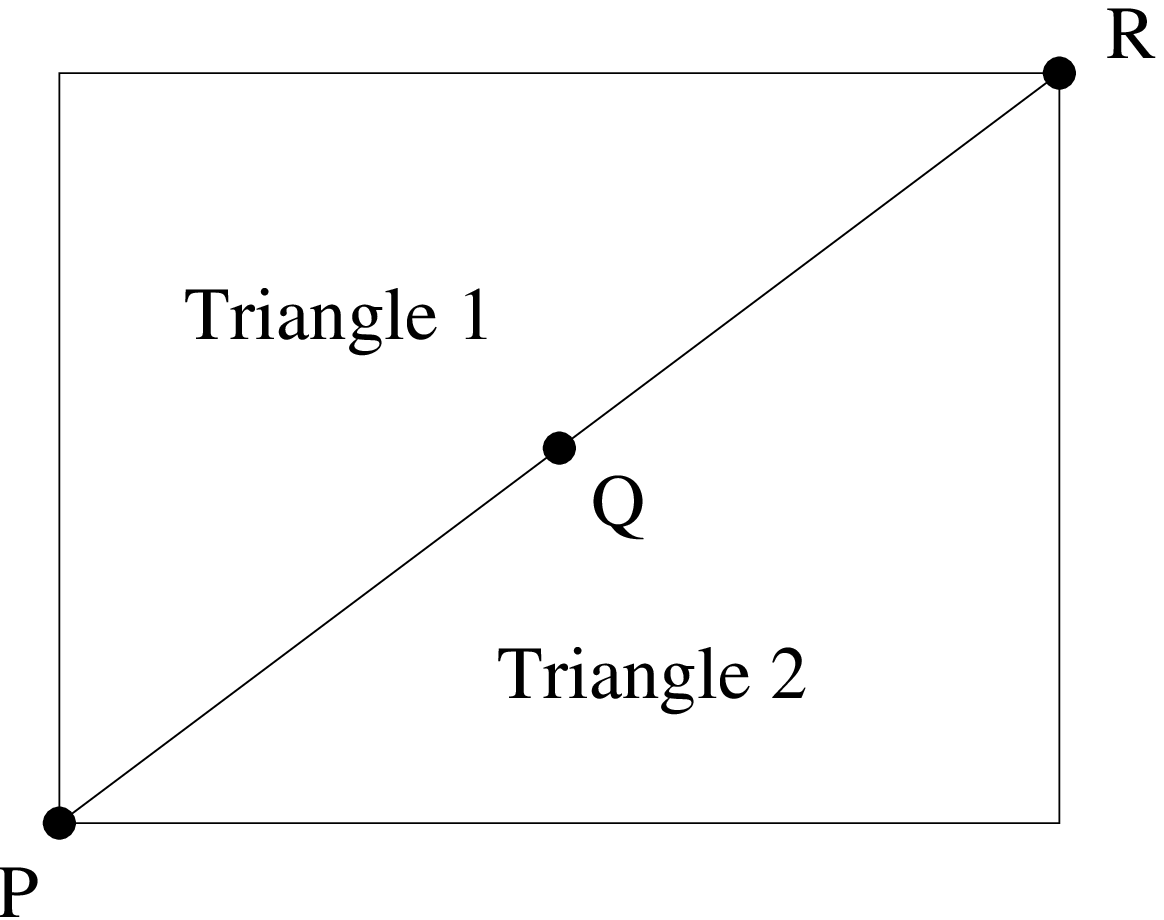}} 
\vs .1
\centerline{Figure 7. A fault line.}
\vs .1

Subdivide the pair of triangles $2n$ times.  The pattern of
subdivision of triangle 2 along the interval $PQ$ is the same as that
of triangle 1 along the interval $RQ$.  We have already seen that only
large triangles abut the main diagonal, and that they do so only along
their hypotenuses and long legs.  For each point $x$ on the interval
$PQ$, let $g_n(x)$ be the number of complete long legs on the triangle
1 side of $Px$ minus the number of complete long legs on the triangle
2 side of $Px$.

{\nd \bf Lemma 8:} {\it $|g_n(Q)| \ge {c \over b} |f(n)| - 1$.  
In particular, $Lim_{n \to \infty} |g_n(Q)| = \infty$.}

\nd {\it Proof:} 
Consider the the $2n$-fold subdivision of the pair of triangles in
Figure 7.  The hypotenuse of triangle 2 gets divided into $\sigma^n(H^+)$.
Let $M_n$ be the point on the hypotenuse corresponding to the middle of
the sequence $\sigma^n(H^+)$.  By construction, there are $2 f(n)$ more
long legs, and $f(n)$ fewer hypotenuses, between $P$ and $M_n$ than between
$M_n$ and $R$.  Since legs have length $b$, hypotenuses have length $c$, and
$2b>c$, $M_n$ lies a distance $f(n)(2b-c)$ closer to $R$ than to $P$, or
a distance $f(n)(2b-c)/2$ beyond $Q$.

Suppose $f(n) < 0$.  Then $M_n$ lies between $P$ and $Q$.  The number of
hypotenuses between $P$ and $Q$ is at least the number of hypotenuses
between $P$ and $M_n$, and so is at least $-f(n)$ more than the number of
hypotenuses between $Q$ and $R$.  Since the length of $PQ$ equals that
of $QR$, and length equals $b \times$ legs plus $c \times$ hypotenuses, 
there are at least $-(c/b) f(n)$ more legs in $QR$ than in $PQ$.  Thus 
$|g_n(Q)|$ is at least the integer part of $(c/b)|f(n)|$, which is greater
than $(c/b)f(n)-1$.  

If $f(n)>0$ this situation is reversed.  Then $M_n$ lies between $Q$
and $R$, and the interval between $M_n$ and $R$ has the surplus of
$f(n)$ hypotenuses.  The interval $QR$ has a surplus of at least that
many hypotenuses, so the interval $PQ$ has a surplus of at least
$(c/b)f(n)$ legs.  Taking the integer part, we see that $PQ$ has a
surplus of at least $(c/b)f(n)-1$ complete legs.  \qed

We are now in a postition to prove Theorem 5.  Let $x$ be a point on
$PQ$ that is the vertex of a tile in the subdivision of triangle 2,
and let $y$ be the nearest vertex, on $PR$, of a tile in the 
subdivision of triangle 1.  The interval $xy$ is the interval of contact
between two tiles, one in triangle 1 and the other in triangle 2.  If 
$y$ lies between $P$ and $x$, then the length of $xy$ is congruent,
modulo $c$, to $b g_n(x)$.  If $x$ lies between $P$ and $y$, then the
length is congruent either to $-b g_n(x)$ or $b-bg_n(x)$ (mod $c$), 
depending on whether $xy$ is part of the hypotenuse or long leg of a tile
on side 1.    

As we pick points $x$ from $P$ to $Q$, $g_n(x)$ goes, by steps of one,
from $0$ to $g_n(Q)$.  Thus the intervals $xy$ take on at least
$|g_n(Q)|/2$ distinct lengths.  Since $|g_n(Q)|$ grows with $n$, all
we need for Theorem 5 is to show that fault lines modeled on
$\sigma^n(H^+)$, for arbitrarily large values of $n$, occur in the
tiling $Til(1/2)$.

By construction, $Til(1/2)$ contains supertiles modeled on $T_m$, for
arbitrarily large values of $m$.  If $m$ is even, then the primary
fault line exhibited in Figure 6 is modeled on $\sigma^n(H^+)$, with
$m = 2n+2$.  If $m$ is odd, then there is also a supertile modeled on
$T_{m-1}$, namely the descendants of the large tile in the first
subdivision of $T$.  Since $m-1$ is even, this supertile contains a
fault line modeled on $\sigma^n(H^+)$, with $m = 2n+3$.  Since $m$ is
unbounded, we have obtained our requisite arbitrarily long fault
lines.

This completes the proof of Theorem 5. \qed

{\nd \bf \S 4. Statistical properties of rational tilings}
 
In this section we consider the statistical distribution of sizes and
orientations of tiles in the various rational tilings $Til(p/q)$.
More precisely, we consider the distribution of sizes and orientations
in a supertile modeled on $T_n$, and take the limit as $n \to \infty$.
We first show that the distribution of sizes approaches a simple limit
as $n \to \infty$.  The limiting distribution is given by the
eigenvector associated to the largest eigenvalue of a $p \times p$ (or
$q \times q$) population matrix.  We also analyze the second
eigenvalue of this matrix.  The failure of the edge-to-edge property
for $Til(1/2)$ was a result of fluctuations that were governed by this
second eigenvalue.  We conjecture that the edge-to-edge property holds
only for those tilings with second eigenvalue smaller than 1, and we
classify these tilings.

We then turn to the joint distribution of sizes and orientations in
rational $p/q$ tilings.  We show that, for each size, the distribution
of orientations is asymptotically uniform.  Specifically, we
parametrize $O(2)$, the group of orientations, by two copies of the
unit circle.  Given an interval in this set, the fraction of tiles in
$T_n$, of a given size, whose orientations lie in that interval,
approaches a constant times the length of the interval.  In the
terminology of Radin [R3], the tilings $T(p/q)$ exhibit ``statistical
rotational symmetry''.

By Theorem 3, $Til(p/q$) contains triangles of $m=\max(p,q)$ distinct
sizes.  $Til(p/q)$ is equivalent to a traditional substitution tiling,
with prototiles $D_1, \ldots D_m$ of $m$ sizes.  We take $D_1$ to be
the largest size and $D_m$ to be the smallest.  Subdivision and linear
rescaling by $r^{-1}$, where $r = (a/c)^{1/p} = (b/2c)^{1/q}$, takes
$D_{i+1}$ to $D_i$, and takes $D_1$ to four copies of $D_q$ and one
copy of $D_p$.  That is, the population matrix, which gives the
population of $T_{n+1}$ in terms of the population of $T_n$, has
matrix elements
$$ M_{ij} = \cases{ 1 & if $j=i+1$; \cr
                    1 & if $j=1$ and $i=p$; \cr
                    4 & if $j=1$ and $i=q$; \cr
                    0 & otherwise.}  \eqno(4.1)
$$

The properties of $M$ are summarized in the following theorem:

\nd {\bf Theorem 6} {\it 
\item {1)} The characteristic polynomial of $M$ is
$$ p(\lambda) = \cases {\lambda^q -  \lambda^{q-p} -4 & if $p<q$; \cr
\lambda^p -  4 \lambda^{p-q} -1 & if $p>q$.} \eqno(4.2) 
$$
\item {2)} The largest eigenvalue of $M$ is $r^{-2}$.
\item {3)} There are exactly $q$ eigenvalues with modulus greater than one.
\item {4)} The eigenvectors $\psi$ of $M$, for fixed eigenvalue $\lambda$, 
take the form
$$ \psi_k = \lambda^k - \lambda^{k-p} H(k-p-1) - 4 \lambda^{k-q}
H(k-q-1),  \eqno(4.3)
$$
where $H(n)$ is the discrete Heavyside function
$$ H(n) = \cases{1 & if $n \ge 0$; \cr 0 & otherwise.} \eqno(4.4) $$
\item {5)} Asymptotically, the number of tiles of  size $D_k$
is a fraction
$$ \nu_k = {1 - r^2 \over 4 c^2} (a^2 H(p-k) + b^2 H(q-k))r^{-2k} \eqno(4.5)
$$
of the total.
\item {6)} Asymptotically, the area covered by tiles of size $D_k$
is a fraction
$$  \rho_k = {a^2 H(p-k) + b^2 H(q-k) \over p a^2 + q b^2} \eqno(4.6)
$$
of the total.}

Note that, if $p<q$, then $H(k-q-1)$ is identically zero and $H(q-k)$
is identically one.  If $p>q$, then $H(k-p-1)$ is identically zero and 
$H(p-k)$ is identically one.  As written, expressions (4.3), (4.5) and 
(4.6) apply to both the $p<q$ and $p>q$ cases.

\nd {\bf Corollary:} {\it Let $t$ be a tile of size $k$.  Let $N_k(n)$
be the number of descendants of $t$ after the $n$-th application
of the substitution rule.  Then 
$$\lim_{n \to \infty} r^{2n} N_k(n) = {4 c^2 r^{2k} 
\over (1-r^2)(pa^2 + qb^2)}.  \eqno(4.7)
$$}

\nd {\it Proof of Corollary:} From the distribution (4.5), 
we compute the average area per tile to be $(1-r^2)ab(pa^2 +qb^2)/8c^2$.
A tile of size $D_k$, subdivided and rescaled $n$ times, has area
$a b r^{2k-2n}/2$.  Dividing by the area per unit tile we obtain (4.7). \qed

\nd {\bf Remark:} If $q>1$, then the second largest eigenvalue of $M$
is greater than one.  The fluctuations 
in population associated to the corresponding
eigenvector then grow with subdivision, although they do not grow as fast as
the population as a whole.  It was precisely this phenomenon that caused
the edge-to-edge property to fail for $Til(1/2)$.  
If $q=1$, then the second eigenvalue is less than one.  In \S 6 we shall
see that, in $Til(2)$, this causes the triangles to meet in 
only a finite number of ways.

\nd {\bf Conjecture:} {\it If $p$ and $q$ are relatively prime integers, $q>1$
and $p/q \ne 1/3$, then the tiling $Til(p/q)$ has tiles that meet in
an infinite number of different ways.  If $p>1$, then the tiles in
$Til(p)$ meet in only a finite number of ways.}
  
\nd {\it Proof of Theorem 6:} For $k<m$, the $k$-th row of the 
vector equation $M \psi = \lambda \psi$ reads
$$ \psi_{k+1} = \lambda \psi_k  - (\delta_{k,p} + 4 \delta_{k,q}) \psi_1.
\eqno(4.8)
$$
Setting $\psi_1=\lambda$, we repeatedly use (4.8) to obtain expression
(4.3) for $\psi_2, \ldots \psi_m$.  Plugging this into the 
$m$-th row of $M \psi = \lambda \psi$ then 
gives the characteristic polynomial (4.2).

Now suppose $p<q$, and consider the function $p(\lambda)$ for $\lambda$ 
real and positive.  Note that $p(r^{-2}) = r^{-2q} - r^{2p-2q} - 4
= 4c^2/b^2 -4a^2/b^2 - 4 =0$, since $a^2+b^2=c^2$.  
When $\lambda \ge 1$, $p'(\lambda) = q\lambda^{q-1}
- (q-p) \lambda^{q-p-1} >0$.  Thus $p(\lambda)>0$ for all $\lambda>r^{-2}$.  

Now let $\lambda \ne r^{-2}$ be a root of $p(\lambda)$.  We will show that
$|\lambda|<r^{-2}$.  If $\lambda$ is real and positive, 
then $\lambda<r^{-2}$.
If $\lambda$ is not real and positive, then, since $q$ and $q-p$ 
are relatively prime, we cannot have $\lambda^q$ and $\lambda^{q-p}$ both
real and positive.  Thus, by the triangle inequality, 
$$ 0 = |p(\lambda)| =  |\lambda^q - \lambda^{q-p} - 4| 
> |\lambda^q| - |\lambda^{q-p}| -4 
=  p(|\lambda|). \eqno(4.9)
$$ 
Since $p(|\lambda|)<0$, we must have $|\lambda| < r^{-2}$.

We count the number of roots in the unit circle via the argument
principle, tracking the argument of $p(\exp(i\theta))$ as $\theta$
goes from $0$ to $2 \pi$.  When $|\lambda|=1$, $p(\lambda)$ always has
negative real part, as $4 > |\lambda^q + \lambda^{q-p}|$.  Thus the
winding is zero, and none of the roots of $p(\lambda)$ lie in the unit
circle.  Thus there are $q$ roots, counted with multiplicity, outside
the unit circle. Since $p'(\lambda)$ is never zero outside the unit
circle, all these roots are distinct.

Now suppose $p>q$.  Then $p(\lambda) = \lambda^{p-q}(\lambda^q -
\lambda^{q-p} -4)$.  By the same arguments as before, $p(r^{-2})=0$,
and $p(\lambda)>0$ for real $\lambda>r^{-2}$.  For any eigenvalue
$\lambda$ other than $r^{-2}$, $0 = |p(\lambda)| > p(|\lambda|)$.
Since $p(|\lambda|)<0$, $|\lambda| < r^{-2}$.  By the argument
principle, there are $p-q$ roots inside the unit circle, since the
dominant term of $p(\lambda)$ on the unit circle is $4 \lambda^{p-q}$.
This leaves $q$ roots outside the unit circle.  This completes the
proof of statements 1--4.

To obtain the asymptotic distribution of sizes, we must decompose the
initial population into eigenvectors of $M$.  The asymptotic distribution
will be the eigenvector corresponding to the largest eigenvalue, assuming
the coefficient of that eigenvector is nonzero.  We have shown that this
largest eigenvalue is $r^{-2}$.  Since the total area
of the system grows by a factor of $r^{-2}$ each time, the coefficient of
this eigenvector is not zero.  Thus the asymptotic distribution of 
population is given by a multiple of the expression (4.3), 
with $\lambda = r^{-2}$.  Normalizing, we obtain
expression (4.6).   Multiplying this by the area of each tile and 
normalizing again gives the asymptotic distribution of areas (4.7).  \qed

We now turn to the joint distribution of sizes and orientation.  To do
this we must first parametrize the space of possible orientations of a
single size.  This space is isomorphic to 2 copies of the unit circle.
We specify both the handedness of the triangle and the direction a
fixed vector in it points in the plane.  We take as our reference
vector the ray from the small angle to the right angle.  In Figure 8,
the first triangle has orientation $(+,\phi_1)$, while the second has
orientation $(-,\phi_2)$.  We will let the $\Omega$ denote the ordered
pair $(\pm, \phi)$, and let $d\Omega = d\phi/4\pi$ be the Haar measure
on the space of orientations.  The space of all possible tiles up to
translation, which we denote $X$, is $2m$ copies of $S^1$.

\vs .05
\hs1.15 \vbox{\epsfxsize=3truein\epsfbox{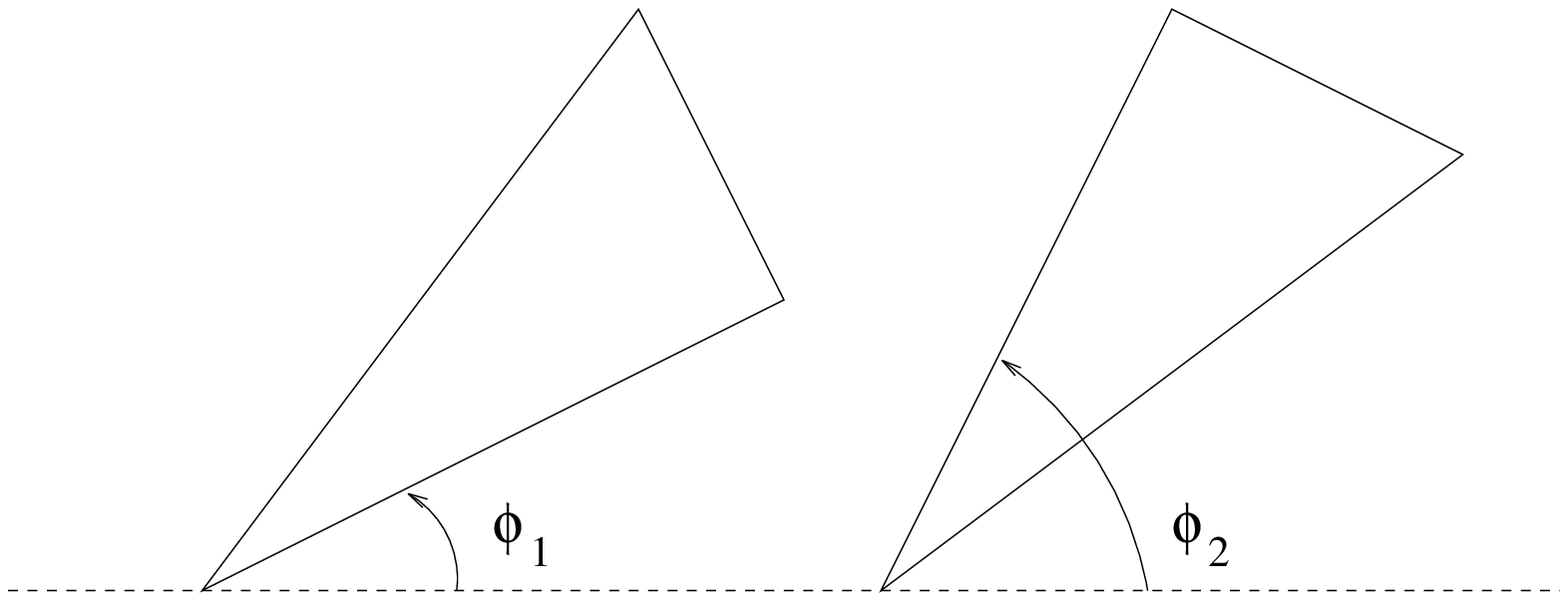}} 
\vs .05
\centerline{Figure 8. Orientation of triangles.}
\vs .1

Let $f$ be a function on $X$.  Given a collection $S$ of tiles, let
$<f,S>$ be the average value of $f$ on the individual tiles in $S$,
where each tile is given equal weight.  Let $<f,S>'$ be the average
value of $f$ on the individual tiles in $S$, where each tile is
weighted proportionally to its area.

\nd {\bf Theorem 7} {\it Assume a tiling $Til(p/q)$, with $p$ and $q$
relatively prime and $p/q \ne 1/3$.
Let $f$ be a continuous function on $X$, and
let $\{ S_n \}$ be a sequence of supertiles of increasing size.  Let
$d \nu = \sum_k \nu_k d\Omega_k$ and $d \rho = \sum_k \rho_k d \Omega_k$
be measures on $X$, where $\nu_k$ and $\rho_k$ are as in Theorem 6.  Then
$$\lim_{n \to \infty} <f,S_n> =  \int_X  f d\nu \eqno(4.10a) $$
$$ \lim_{n \to \infty} <f,S_n>' =  \int_X  f d\rho. \eqno(4.10b)
$$
These limits also apply if $f$ is the characteristic function of an
interval in $X$.}

\nd {\it Proof:} We first reduce the problem to establishing (4.10a)
for an arbitrary continuous function $f$.  Once we have established
(4.10a), (4.10b) follows by applying (4.10a) to the continuous
function $\tilde f(k,\Omega) = r^{2k} f(k,\Omega)$.  Once equations
(4.10ab) have been established for continuous functions, the extension
to characteristic functions is standard (for details see [CFS]).  Let
$I$ be an interval, and let $\chi_I$ be its characteristic function.
We choose continuous functions $f^{\pm}_\epsilon$, such that
$f^-_\epsilon \le \chi_I \le f^+_\epsilon$, and such that
$\lim_{\epsilon \to 0} \int_X f^{\pm}_\epsilon d \nu = \int_I d\nu$.
$\lim_{n \to \infty} <\chi_I,S_n>$ is sandwiched between $\lim_{n \to
\infty} <f^-_\epsilon,S_n>$ and $\lim_{n \to \infty}
<f^+_\epsilon,S_n>$, hence between $\int_X f_\epsilon^- d\nu$ and
$\int_X f_\epsilon^+ d\nu$, and so must equal $\int_I d\nu$.

To establish (4.10a) we must introduce some notation.  If $f$ is a
function on $X$ and $S$ is a collection of tiles, let $(f,S)$ be the
sum of $f$ evaluated on the individual tiles of $S$.  Let $\Phi$
denote the action of subdividing and rescaling.  That is, $\Phi$
acting on a tile of size $D_{k+1}$ gives a tile of size $D_k$ of the
same orientation, while $\Phi$ acting on a tile of size $D_1$ gives
one tile of size $D_p$ and four of size $D_q$, having various
orientations.  Let $\Phi^*$ be the dual of $\Phi$ by $(\cdot,\cdot)$,
acting on the space of functions:
$$ (\Phi^*f,S) = (f,\Phi(S)). \eqno(4.11) $$
Let $f_n = (\Phi^*)^n f$.
Note that $\Phi^*$ is linear and sends non-negative functions to 
non-negative functions, so if $f \le g$, then $f_n \le g_n$.  

Now suppose we have a sequence of tiles $t_n$ and supertiles 
$S_n = \Phi^n t_n$.  We have that
$$ <f,S_n> = {(f,S_n) \over \hbox{\# of tiles in $S_n$}} = 
{(f_n,t_n) \over \hbox{\# of tiles in $S_n$}}. \eqno(4.12)
$$
Since the number of tiles is given asymptotically by (4.7),
(4.10a) is equivalent to 
$r^{2n} f_n$ converging uniformly to 
$\int_X  f d\nu$ times
$$ \zeta_0(k,\Omega) = {4 c^2 r^{2k}\over (1-r^2)(pa^2 + qb^2)}. 
\eqno(4.13) $$ 

We examine the spectrum of the linear operator $\Phi^*$ 
on the function space $C(X)$. 
The key lemma, whose proof we defer, is

\nd {\bf Lemma 8} {\it The spectrum of $\Phi^*$ is pure point. $\zeta_0$
is an eigenfunction with eigenvalue $r^{-2}$.  All other eigenvalues
have norm strictly less than $r^{-2}$.  Any continuous function $f$
can be written as a (possibly infinite) sum of eigenfunctions of
$(\Phi^*)^{|q-p|}$, such that a subsequence of partial sums converges
uniformly to $f$.}

\nd {\bf Remark:} If $p<q$, it turns out that there are a number of 
functions $\zeta$ for which $(\Phi^*)^{|q-p|}
\zeta = 0$, but $\Phi^*\zeta \ne 0$.  Thus, to achieve a
basis for the 
space of continuous functions, we must use eigenfunctions of $(\Phi^*)^{|q-p|}$
rather than just eigenfunctions of $\Phi^*$.
 
Given the lemma, we write 
$$ f = \sum_{i=0}^\infty c_i \zeta_i. \eqno(4.14) $$
Since a subsequence of the partial sums converges uniformly, and since
$\zeta_0$ has a positive minimum, for each $\epsilon>0$ we can find an
integer $N$ such that
$$ 
f^- \equiv \left ( \sum_{i=0}^N c_k \zeta_i \right ) - \epsilon \zeta_0
< f <  \left ( \sum_{i=0}^N c_k \zeta_i \right ) + \epsilon \zeta_0
\equiv f^+, \eqno(4.15) 
$$
where each $\zeta_i$ is an eigenfunction with eigenvalue $\lambda_i$.
(Strictly speaking, $\zeta_i$ is merely an eigenfunction of
$(\Phi^*)^{|q-p|}$, not necessarily of $\Phi^*$, but this distinction
makes no difference).  Since $|r^{-2}\lambda_i|<1$ for all $i>0$, for
$n>|p-q|$ we have
$$ (c_0 - \epsilon) \zeta_0 + \sum_{i=1}^N c_k (r^2 \lambda_i)^n \zeta_i
\le r^{2n} f_n \le (c_0 + \epsilon) \zeta_0 + \sum_{i=1}^N c_k (r^2 \lambda_i)^n \zeta_i. \eqno(4.16)
$$
As $n \to \infty$, the left hand side converges uniformly to $(c_0-\epsilon)
\zeta_0$, while the right hand side converges uniformly to $(c_0+\epsilon)
\zeta_0$.  Since $\epsilon$ is arbitrary, $r^{2n}f_n$ must converge
uniformly to $c_0 \zeta_0$.

All that remains is to compute $c_0$ in terms of $f$. Since $\zeta_0$
is invariant under rotation and reflection, $c_0$ must be of the form
$\sum_k d_k \int f(k,\Omega) d\Omega$ for some universal constants
$d_k$.  By comparing characteristic functions of different sizes, we
see that the constants $d_k$ must be proportional to $\nu_k$.
Finally, for the constant function $f=1$, $f_n(k,\Omega)$ is the
number of descendants of a tile of size $D_k$, which we have already
computed in (4.7). This fixes the proportionality constant. \qed

\nd {\it Proof of Lemma 8:}  $\Phi^*$ commutes with rotations, so we may
simultaneously diagonalize $\Phi^*$ and the rotation operator $-i {d
\over d\phi}$.  The eigenvalues of the rotation operator are of course
the integers, with each eigenspace being $2m$-dimensional.
Specifically, the eigenspace corresponding to an integer $n$ is the
span of the $2m$ functions obtained by restricting $\exp(i n \phi)$ to
each of the $2m$ circles in $X$. Operators on finite dimensional
spaces always have pure point spectra.  Summing over $n$, we see that
the spectrum of $\Phi^*$ is pure point.

On each $2m$ dimensional subspace corresponding to the Fourier mode
$\exp(i n \phi)$, the action of $\Phi^*$ is described by a $2m$ by
$2m$ matrix $E$, which we write as an $m \times m$ array of $2 \times
2$ matrices.  Let $\theta = \tan^{-1}(a/c)$ be the acute angle in our
basic triangle.  Consider the matrices
$$ A = \left ( \matrix{0 & \exp[-in(\theta  + \pi/2)] \cr
\exp[in(\theta  + \pi/2)] & 0} \right ), \eqno(4.17)
$$
$$
B = \left ( \matrix{\exp(in\theta ) + \exp(i n (\theta +\pi)) 
& 2 \exp(-in\theta ) \cr
2 \exp(in\theta ) &  \exp(-in\theta ) + \exp(i n (-\theta +\pi))} \right ).
\eqno(4.18) $$
$A$ and $B$ describe the orientations of the five daughter tiles in
terms of the orientation of the parent tile, as expressed in the
$n$-th representation of the rotation group $SO(2)$.  Specifically,
$A$ describes the daughter tile of hypotenuse $a$, while $B$ describes
the four daughter tiles of hypotenuse $b/2$.  In each case the first
column describes the daughters of a positively oriented tile, while
the second column describes the daughters of a negatively oriented
tile.  See Figure 9.

\vs .1
\hs1.15 \vbox{\epsfxsize=3truein\epsfbox{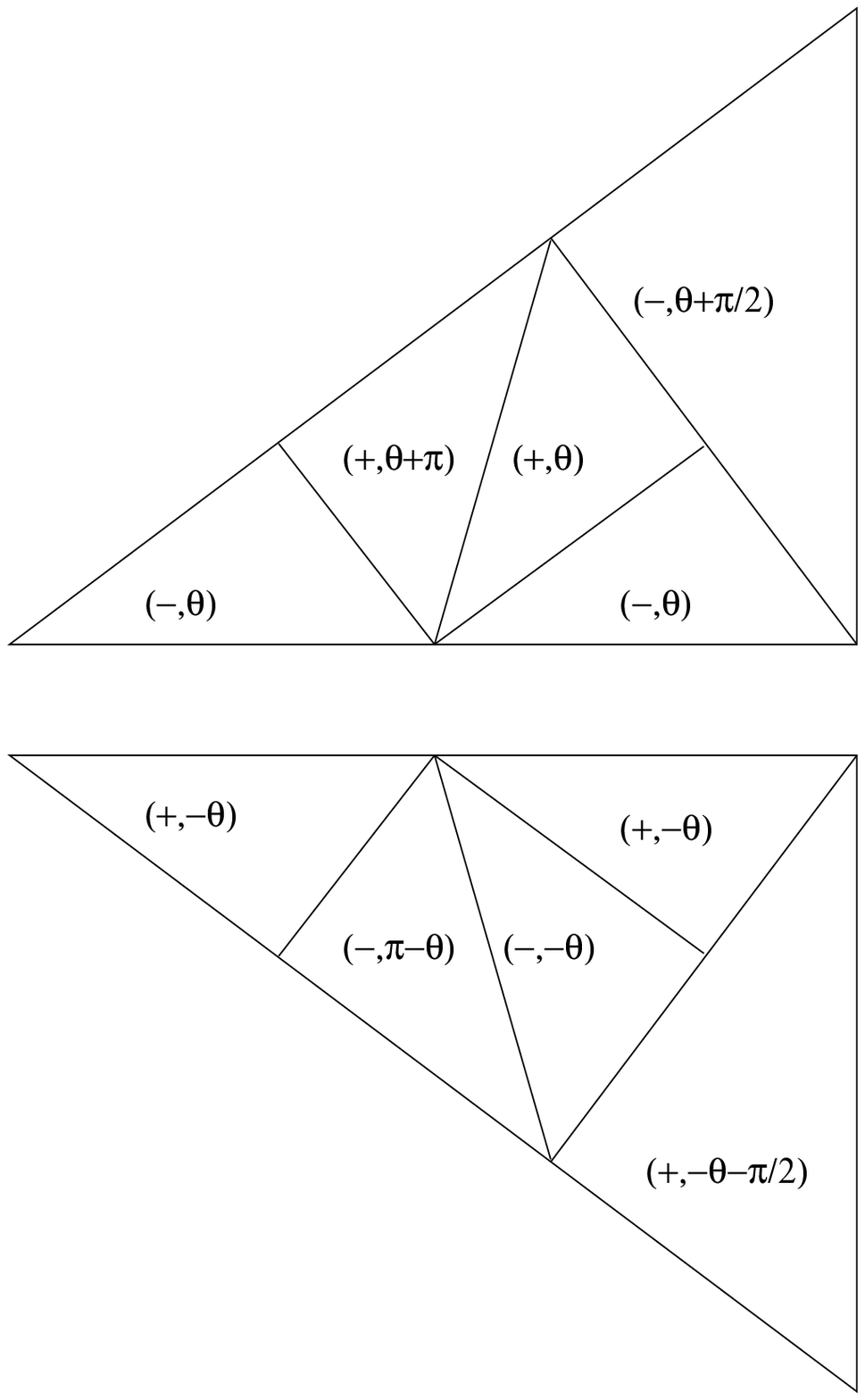}} 
\vs .1
\centerline{Figure 9. Orientation of daughter tiles.}
\vs .1

Our matrix $E$ has the following matrix elements:
$$ E_{ij} = \cases{ 1 & if $j=i-1$; \cr
                    A & if $i=1$ and $j=p$; \cr
                    B & if $i=1$ and $j=q$; \cr
                    0 & otherwise.}  \eqno(4.19)
$$
Notice that $E$ is essentially the transpose of our population matrix
$M$, with the daughters represented by the matrices $A$ and $B$ rather
than the numbers 1 and 4.  We compute the eigenvalues and eigenvectors
of $E$ by directly solving the equation $E \psi = \lambda \psi$.

Suppose $\psi_m = v$, where $v$ is a 2 component vector. Each row but
the first of $E \psi = \lambda \psi$ implies that $\psi_{k-1} =
\lambda \psi_k$, and hence that $\psi_k = \lambda^{m-k} v$ for all
$k$.  The first row then says that
$$ (\lambda^m - \lambda^{m-q} B - \lambda^{m-p} A)v= 0. \eqno(4.20)
$$
Taking the determinant of the matrix on the left hand side, we obtain
the characteristic polynomial of $E$,
$$ \eqalign{ p(\lambda) =  & \lambda^{2m} - \lambda^{2m-q} (2\cos(n\theta ) + 2\cos(n\pi+n\theta ))
 + \lambda^{2m-2q}(2\cos(n\pi)-2) 
\cr & - \lambda^{2m-2p} - \lambda^{2m-p-q}
(4\cos(n \pi/2)).} \eqno(4.21)
$$

First we consider $n=0$.  If $p<q$, $p(\lambda)= (\lambda^q -
\lambda^{q-p} -4) (\lambda^p+1) \lambda^{q-p}$.  We recognize the
first factor as the characteristic polynomial of $M$, with one root
$r^{-2}$ and all other roots smaller in norm.  The roots of the second
and third factors have norm 1 and 0, respectively.  The eigenvector
corresponding to eigenvalue $r^{-2}$ is $\psi_k = r^{2k}\left
(\matrix{1 \cr 1} \right )$, i.e., $\zeta_0$ (up to scale).

We can understand the eigenvalues as follows.  The roots of the first
factor all correspond to $v = \left (\matrix{1 \cr 1} \right )$, and
describe fluctuations in the numbers of tiles of various sizes,
irrespective of orientation.  This is the problem we previously
studied in Theorem 6.  The unit and zero eigenvalues correspond to $v
= \left (\matrix{1 \cr -1} \right )$.  The zero eigenvalue has
algebraic multiplicity $q-p$ and geometric multiplicity 1.  Since the
basic subdivision produces equal numbers of positively and negatively
oriented ``B'' tiles, after $q-p$ subdivisions there must be the same
number of positively and negatively oriented tiles of each size
$D_{p+1}, \ldots D_m$.  Since there is but one ``A'' daughter of each
subdivided tile, the imbalance between positive and negative
orientation in the larger sizes neither grows nor shrinks, but
oscillates with period $2p$. These correspond to the roots of
$\lambda^p+1$.

If $p>q$, then $p(\lambda)$ factorizes as $p(\lambda)=(\lambda^p - 1
-4\lambda^{p-q})(\lambda^p+1)$.  Again, the first factor is the
characteristic polynomial of $M$, and governs the total number of
tiles of each size, with all eigenvectors having $v = \left (\matrix{1
\cr 1} \right )$.  The largest eigenvalue is $r^{-2}$, with
eigenvector $\zeta_0$.  The roots of $(\lambda^p+1)$ describe
oscillations in the numbers of positively vs. negatively oriented
tiles and have $v = \left (\matrix{1 \cr -1} \right )$.
 
Next we consider $n$ odd.  Then $\exp(in\pi)=-1$, and $p(\lambda)$
simplifies to $\lambda^{2m} - 4 \lambda^{2m-2q} - \lambda^{2m-2p}$.
This is just the characteristic polynomial of $M$ applied to
$\lambda^2$.  By Theorem 6, the largest roots have $\lambda^2=r^{-2}$,
or $\lambda= \pm r^{-1}$.

Finally we consider $n$ even and nonzero.  Then $p(\lambda)= 
\lambda^{2m} - \lambda^{2m-q}(4 \cos(n \theta )) - \lambda^{2m-2p} - 
\lambda^{2m-p-q} (4 \cos(n \pi/2))$.  We show that all roots are
smaller than $r^{-2}$ by
the argument principle.  Note that $\theta $ is an irrational multiple
of $\pi$, so $\cos(n\theta ) \ne \pm 1$. On the circle
$|\lambda|=r^{-2}$ we have
$$ 
\eqalign{ |\lambda^{2m-q} (4 \cos(n \theta ))| + |\lambda^{2m-2p}| + 
& |\lambda^{2m-p-q}(4 \cos(n \pi/2))| \cr =  & |\lambda^{2m}| \left (
{b^2 |\cos(n \theta )| \over c^2} + {a^4 \over c^4} + {a^2 b^2 \over
c^4}
\right ) \cr
= & |\lambda^{2m}| {a^2 + b^2 |\cos(n \theta )| \over c^2} \cr < &
|\lambda^{2m}|.} \eqno(4.22) $$ Thus the $\lambda^{2m}$ term of
$p(\lambda)$ dominates on the circle $|\lambda|=r^{-2}$, so the
winding of the argument of $p(r^{-2}\exp(i\a))$, as $\a$ goes from $0$
to $2\pi$, is $2m$, and all $2m$ eigenvalues of $E$ lie inside the
circle of radius $r^{-2}$.

This completes the analysis of the spectrum of $\Phi^*$.  Now we need
only consider the decomposition into eigenvectors.  Since $f$ is
continuous, it has an absolutely summable Fourier series.  So we may
write $f(k,\Omega) = \sum_{n} c_n \psi^{(n)}_k \exp(i n \phi)$, where
the $c_n$'s are absolutely summable and $\psi_k^{(n)}$ is a vector in
$\C^{2m}$ with no component larger than 1.  In particular, the partial
sums converge uniformly to $f$.  But $\Psi_k$ is itself a sum of
eigenvectors of $E^{|q-p|}$, so we may rewrite our sum as a sum of
eigenfunctions of $(\Phi^*)^{|q-p|}$.  \qed

{\nd \bf \S 5. Irrational tilings}
 
In this section we consider tilings $Til(z)$, with $z$ irrational.
The analysis is formally similar to that of \S 4, except that we are
now dealing with an infinite number of possible sizes.  In place of
the discrete size parameter $k$ we introduce a continuous size
parameter $s$.  In place of the discrete evolution operator $\Phi^n$
we introduce a 1-parameter semigroup $e^{tL}$.  Although the
continuous case is technically more difficult than the discrete cases,
the results are extremely similar.  Indeed, if one has a sequence of
rational numbers $p_i/q_i$ converging to the irrational number $z$,
then the statistical properties of $Til(z)$ may be obtained as limits
of the corresponding properties of $Til(p_i/q_i)$.  (Note that the
reverse does not hold.  One cannot obtain the statistics of a rational
tiling by taking a limit of irrational tilings.)

As always, we consider a basic right triangle $T_0$ with sides $a$ and
$b$ and hypotenuse $c$.  Let $\a = \ln(c/a)$, $\b=\ln(2c/b)$, and
assume that $z=\a/\b$ is irrational.  By a triangle of size $s$, we
mean a triangle, similar to $T_0$, with hypotenuse $c e^{-s}$.  Note
that larger values of $s$ correspond to smaller triangles, just as in
the rational case, where the size $D_k$ of triangles decreased with
$k$.  In our tiling the size parameter $s$ will take values in
$[0,\mu)$, where $\mu=\max(\a,\b)$.

We now describe a semigroup similar to $\Phi^n$. Let $S$ be a
collection of tiles, all with size in $[0,\mu)$.  Expand this
collection by a linear factor $e^t$, resulting in triangles with sizes
in $[-t,\mu-t)$.  Then subdivide the largest triangle, subdivide the
largest remaining triangle, and so on, until all triangles have
non-negative size parameter.  By Lemmas 1 and 2, this occurs in a
finite number of steps, and results in a collection of tiles with
sizes in $[0,\mu)$. This collection is $e^{tL}S$.

The semigroup $e^{tL}$ naturally acts on the distribution of sizes.
The properties of this action are summarized by the following theorem,
which should be compared with Theorem 6:

\nd {\bf Theorem 8} {\it 
\item {1)} The eigenvalues $\lambda$ of $L$ are the roots of
$p(\lambda)=0$, where
$$ p(\lambda) = e^{\mu\lambda} - e^{(\mu-\a)\lambda} - 4
e^{(\mu-\b)\lambda}. \eqno(5.1) 
$$
There are no multiple eigenvalues.
\item {2)} $\lambda=2$ is an eigenvalue.  All other eigenvalues have
real part strictly less than 2.  If $\a<\b$, all eigenvalues have
real part greater or equal to the real root of 
$e^{\b\lambda} + e^{(\b-\a)\lambda} - 4$, 
while if $\a>\b$, all eigenvalues have
real part greater or equal to the real root of 
$e^{\a\lambda} + 4
e^{(\a-\b)\lambda} - 1$.
\item {3)} The eigenfunction $\psi(s)$, for fixed eigenvalue $\lambda$, 
take the form
$$ \psi(s) = e^{\lambda s} - e^{\lambda (s-\a)} h(s-\a) - 4 e^{\lambda (s-\b)} h(s-\b),  \eqno(5.2)
$$
where $h(x)$ is the Heavyside function
$$ h(x) = \cases{1 & if $x \ge 0$; \cr 0 & otherwise.} \eqno(5.3) $$
\item {4)} Given an interval $I \subset [0,\mu)$, the number of tiles 
with size in $I$ is asymptotically a fraction
$$ {1 \over 2 c^2} \int_I ds (a^2 h(\a-s) + b^2 h(\b-s))e^{2s} 
\eqno(5.4)
$$
of the total.
\item {5)} Given an interval $I \subset [0,\mu)$, the area covered by tiles 
with size in $I$ is asymptotically a fraction 
$$ {1 \over a^2 \a + b^2 \b} \int_I ds a^2 h(\a-s) + b^2 h(\b-s) 
 \eqno(5.5)
$$
of the total.}

\nd {\bf Corollary:} {\it Let $T$ be a tile of size $s$.  Let $N_t(s)$
be the number of tiles in $e^{tL}(T)$.  Then 
$$\lim_{t \to \infty} e^{-2t} N_t(s) = {2 c^2 e^{-2s} 
\over a^2 \a + b^2 \b}.  \eqno(5.6)
$$}

\nd {\it Proof of Corollary:} From the distribution (5.4), 
we compute the average area per tile to be $ab(a^2\a+b^2\b)/4c^2$.
A tile of size $s$, rescaled by a factor $e^t$, has area
$a b e^{2t-2s}/2$.  Dividing by the area per unit tile we obtain (5.6). \qed
  
\nd {\it Proof of Theorem 8:} For $0<t<\min(\a,\b,|\a-\b|)$, the action
of $e^{tL}$ on population distribution functions is
$$ \left( e^{tL} \psi\right) (s) = \cases{
\psi(s+t) + 4 \psi(s+t-\b) & if $s \in [\b-t,\b)$; \cr
\psi(s+t) + \psi(s+t-\a) & if $s \in [\a-t,\a)$; \cr
\psi(s+t) & all other $s \in [0,\mu)$,} 
\eqno(5.7) $$
where $\psi(s)$ and $e^{tL}\psi(s)$ are understood to be zero 
for $s \ge \mu$ or $s<0$.  $e^{tL}$ acts continuously on
functions that are continuous away from $0,\a,\b$, and for which
the oscillations at $\a$ and $\b$ are given by
$$ \psi(\a^-) - \psi(\a^+) = f(0); \qquad \psi(\b^-) - \psi(\b^+) = 4 f(0).
\eqno(5.8)
$$

Setting $e^{tL}\psi = e^{t\lambda}\psi$ we see that the eigenfunction
$\psi(s)$ must equal $e^{\lambda s}$ times a piecewise constant
function with discontinuities at $0,\a,\b$.  Applying the boundary
conditions (5.8), we obtain the eigenfunction (5.2).  For $s>\mu$,
$\psi(s)$ then equals $e^{(s-\mu)\lambda}$ times $p(\lambda)$.  The
vanishing of $\psi(s)$ for $s>\mu$ is equivalent to the eigenvalue
equation $p(\lambda)=0$. Thus eigenfunctions satisfying the boundary
conditions are in 1--1 correspondence with roots of $p(\lambda)$, with
the eigenfunctions given by (5.2).

Suppose $\a < \b$, in which case $p$ takes the form
$$ p(\lambda) = e^{\b\l} - e^{(\b-\a)\l} -4. \eqno(5.9) $$
$\l=2$ is a root, since $p(2)= e^{2\b} - e^{2(\b-\a)}-4 = {4 c^2 \over b^2}
- {4 a^2 \over b^2} - 4 = 0$.  This is the only real root, insofar as $p(\l)$
is an increasing function of $\l$ for $\l>0$, and $p(\l)$ is negative for
$\l \le 0$.  In particular, $p(\l)<0$ implies that $\l < 2$.

Now consider complex roots $\l=\l_R + i \l_I$.  If $\l_I \ne 0$, 
$e^{\b\l}$ and $e^{(\b-\a)\l}$ cannot both be real, insofar as $\b$ is not 
a rational multiple of $\b-\a$.  Thus $|e^{\b\l}|$, $|e^{(\b-\a)\l}|$ and 
4 satisfy a strict triangle inequality.  In particular,
$$ 0 > |e^{\b\l}| - |e^{(\b-\a)\l}| -4 = p(\l_R), \eqno(5.10) $$ 
so $\l_R <2$.  Also,
$$0 <  |e^{\b\l}| + |e^{(\b-\a)\l}| -4 = e^{\b\l_R} + e^{(\b-\a)\l_R} -4,
\eqno(5.11) $$
so $\l_R$ is greater than the real root of $e^{\b\l} + e^{(\b-\a)\l}-4$.

Now we exclude the possibility of multiple roots.  A multiple root
would require $p(\l)=p'(\l)=0$.  Suppose $0=p'(\l)= \b e^{\b\l}
+(\b-\a)e^{(\b-\a)\l}$. Then $e^{\b\l}$ and $e^{(\b-\a)\l}$ must have
the same phase, and their difference must also have that phase.
However, if $\l$ is not real, $e^{\b\l}$ and $e^{(\b-\a)\l}$ cannot
both be real, so their difference is not real, so their difference is
not 4.  Thus $p'(\l)=0$ implies that $p(\l) \ne 0$, and there are no
multiple roots away from the real axis.  On the real axis, the only
root is $\l=2$, and we have already seen that $p'(2)$ is positive, not
zero.

This establishes statements 1--3 for the case $\a < \b$.  The argument
for $\a > \b$ is almost identical, and is not repeated.

Now suppose that we initially have a population distribution that
is a linear combination of the eigenfunctions $\psi_\l$.  Applying
$e^{tL}$ to the system, the total area grows as $e^{2t}$, so the
number of tiles is bounded, both above and below, by a multiple of $e^{2t}$.  
Applying $e^{tL}$ and dividing by the number of tiles damps out all
the modes with eigenvalue less than 2, i.e. all eigenvectors other than
$\psi_2$.  In this case, the final distribution of sizes approaches a 
multiple of $\psi_2$, and statements 4 and 5 of the theorem follow.

Unfortunately, we cannot a priori assume that the initial condition is 
a linear combination of eigenfunctions of $L$, or that a test function is a
linear combination of eigenfunctions of the dual operator $L^*$.  
$e^{tL}$ and $e^{t{L^*}}$ are neither finite-rank operators nor 
self-adjoint operators on a Hilbert space, so standard
theorems about the completeness of a basis is eigenfunctions cannot
be applied.   In principle it is possible for a test function $f$
to have the property that $e^{t{(L^*-2)}} f$ does not converge at all.
We must show that, when $f=\chi_I$  (the characteristic function of an 
interval $I$), $e^{t{(L^*-2)}}f$ {\it does} converge. 

\nd {\bf Lemma 9:} {\it Given an interval $I\in [0,\mu)$.  The fraction of 
the area of $e^{tL}(T_0)$ covered by tiles with size in $I$ approaches a 
limit as $t \to \infty$.}

Given this lemma, it follows that the distribution of area of
$e^{tL}(S)$, for any collection of tiles $S$, approaches a limit, from
which it follows that the distribution of population of $e^{tL}(S)$
also approaches a limit.  Since these limits are invariant under
further evolution, and since the total area is proportional to
$e^{2t}$, these limits must correspond to the $\l=2$ eigenvector of
$L$, hence must take the form (5.4) and (5.5). \qed

\nd {\it Proof of Lemma 9:}  Assume $\a < \b$; the other case is similar.
Let $F_I(t)$ be  the fraction of area of $e^{tL}(T_0)$ 
covered by tiles with size in $I$.  Given an $\epsilon>0$, 
we will show how to compute a number such that, for
all $t$ sufficiently large, $F_I(t)$ is within $\epsilon$
of this number.  Since this can be done for any $\epsilon$, 
$\lim_{t \to \infty} F_I(t)$ must exist.  

It suffices to show that the eventual fraction in an interval of size
$\Delta$, entirely in $[0,\a)$ or in $[\a,\b)$, and with $\Delta$
sufficiently small, can be estimated to 
within $O(\Delta^2)$.  Any larger interval can be broken up into a finite 
number of such small pieces, such that $\sum$ errors $ < \epsilon$.
So let us fix an interval $I$, centered at $s_0$, with width $\Delta$.

The strategy is this:  
We begin with an exact expression for the fraction of area of $e^{tL}(T_0)$
represented by tiles of size $s$.  We sum this over $s \in I$ to get an
exact formula for $F_I(t)$.   By taking
certain limits and replacing certain sums with integrals, we obtain an 
expression that is independent of $t$.  In the process we introduce two
types of errors.  One type can be made arbitrarily small by requiring
$t$ to be sufficiently large.  The other
type is $O(\Delta^2)$.
 
How many triangles of size $s$ appear in $e^{tL}(T_0)$?  That depends on 
whether $s+t$ can be written as $n_1 \a + n_2 \b$ for non-negative $n_1$
and $n_2$.  If $s+t=n_1 \a+n_2 \b$, then a triangle of size $s$
may be obtained by taking a triangle of size $-t$, subdividing it, picking
a daughter, subdividing it, picking a daughter, and so on for $n_1+n_2$
subdivisions, with the descent involving $n_1$ daughters of 
type $A$ and $n_2$ daughters of type $B$.  If
$s\ge \a$, then the last daughter must be of type $B$, or else after
$n_1+n_2-1$ steps we would have already obtained a tile of size in $[0,\b)$,
and would not have made the final subdivision.  If $s<\a$ there is no such
constraint.  We thus have
$$ \hbox{Number of tiles of size $s$} = \cases{
{n_1 + n_2 \choose n_1} 4^{n_2} & if $s+t = n_1 \a + n_2 \b$ and 
$s \in [0,\a)$; \cr
{n_1 + n_2 -1 \choose n_1} 4^{n_2} & if $s+t = n_1 \a + n_2 \b$ and 
$s \in [\a,\b)$; \cr 0 & otherwise.} \eqno(5.12)
$$
Since at each division a fraction $a^2/c^2$ of the area goes into the
$A$ daughter, while a fraction $b^2/c^2$ goes into the $B$ daughters,
then the fraction of the total area represented by tiles of size $s$
is
$$ F_s(t) = \cases{
{n_1 + n_2 \choose n_1} \left ({a^2 \over c^2} \right )^{n_1} 
\left ({b^2 \over c^2} \right )^{n_2} & if $s+t = n_1 \a + n_2 \b$ and $s \in [0,\a)$; \cr
{n_1 + n_2 -1 \choose n_1} \left ({a^2 \over c^2} \right )^{n_1} 
\left ({b^2 \over c^2} \right )^{n_2}  & if $s+t = n_1 \a + n_2 \b$ and $s \in [\a,\b)$; \cr
0 & otherwise.} \eqno(5.13)
$$
Note that $n_2 = (s + t - \a n_1)/\b$.  Now let 
$$f(s,t,n)= \cases{
{(s+t+n(\b-\a))/\b  \choose n} \left ({a^2 \over c^2} \right )^{n} 
\left ({b^2 \over c^2} \right )^{(s + t - \a n)/\b} & if $s \in [0,\a)$; \cr  
{(s+t+n(\b-\a)-\b)/\b  \choose n_1} \left ({a^2 \over c^2} \right )^{n} 
\left ({b^2 \over c^2} \right )^{(s + t - \a n)/\b} & if $s \in [\a,\b)$,}
\eqno(5.14)
$$
and let $\delta_p$ be the periodic $\delta$-function
$$ \delta_p(x) = \sum_{n \in \Z} \delta(x-n). \eqno(5.15) $$
Note that $f(s,t,n)$ is well-defined even when $(s+t+n(\b-\a))/\b$ is not
an integer.  For $t$ large, $f(s,t,n)$ is a slowly-varying function of 
$s$ and $n$.  

We then compute, exactly,
$$F_I(t) = \sum{s \in I} F_s(t) = \sum_{n=0}^{[t/\a]} \int_I
f(s,t,n) \delta_p((s+t-n\a)/\b) \eqno(5.16) 
$$
Next we approximate, by
replacing $f(s,t,n)$ by $f(s_0,t,n)$, where $s_0$ is the midpoint of $I$.
This introduces an error that is a fraction $O(\Delta)$ of the total.
Since the total will turn out to be $O(\Delta)$, the error introduced is
$O(\Delta^2)$.  We thus have
$$F_I(t) = \sum_{n=0}^{[t/\a]} f(s_0,t,n) \int_I \delta_p((s+t-n\a)/\b) + 
O(\Delta^2). \eqno(5.17)
$$

Next we use the fact that multiples of an irrational number are
uniformly distributed on $\R/\Z$.  For any desired degree of accuracy,
one can find an $N$ such that, for any $N$ consecutive integers $n_i$,
$\int_I \delta_p((s+t-n_i\a)/\b)$ equals one a fraction $\Delta/\b$ of
the time (to within the allowed error), and 0 the rest of the time.
If we restrict ourselves to $t$ so large that $f(s_0,t,n)$ is nearly
constant as $n$ varies over intervals of size $N$, then
$$ \sum_{n=0}^{[t/\a]} f(s_0,t,n) \int_I \delta_p((s+t-n\a)/\b)
\approx (\Delta/\b) \sum_{n=0}^{[t/\a]} f(s_0,t,n). \eqno(5.18)
$$
Finally, the $t \to \infty$ limit of the resulting sum can be expressed
as an integral, and yields a nonzero number.          \qed

We now turn, as in \S 4, to the joint distribution of sizes and
orientation.  We parametrize the space of possible orientations of a
single size, as in \S 4, by two copies of the unit circle, with Haar
measure $d\Omega$. (See Figure 8).  The space of all possible tiles up
to translation, which we denote $X$, is $2$ copies of $S^1 \times
[0,\mu)$.

Let $f$ be a function on $X$.  Given a collection $S$ of tiles, let
$<f,S>$ be the average value of $f$ on the individual tiles in $S$,
where each tile is given equal weight.  Let $<f,S>'$ be the average
value of $f$ on the individual tiles in $S$, where each tile is
weighted proportionally to its area.  Let $(f,S)$ be the sum of $f$ on
the individual tiles of $S$.  We define a semigroup $e^{t{L^*}}$
acting on functions on $X$ by
$$ (e^{t{L^*}}f,S) = (f,e^{tL}(S)). \eqno(5.19) $$

\nd {\bf Theorem 9} {\it Assume a tiling $Til(z)$, with $z$ irrational
and with $\theta=\tan^{-1}(a/b)$ an irrational multiple of $\pi$.  Let
$f$ be a continuous function on $X$, and let $\{ S_n \}$ be a sequence
of supertiles of increasing size.  Let $d \nu = e^{2s}(a^2 h(\a-s) +
b^2 h(\b-s)) d\Omega ds /2c^2$ and $d \rho = 2 c^2 e^{-2s} d\nu /
(a^2\a+b^2\b)$ be measures on $X$, where $h$ is the Heavyside function
(5.3).  Then
$$\lim_{n \to \infty} <f,S_n> =  \int_X  f d\nu \eqno(5.20a) $$
$$ \lim_{n \to \infty} <f,S_n>' =  \int_X  f d\rho. \eqno(5.20b)
$$
These limits also apply if $f$ is the characteristic function of a
rectangle in $X$.}

\nd {\bf Remark:} The measures $d\nu$ and $d\rho$ are closely related to the
integrands in (5.4) and (5.5), respectively.  Theorem 9 states that
the joint distribution of sizes and orientations is a product: the
size distribution previously found in Theorem 8 times a uniform
distribution of orientations.

\nd {\it Proof:} The proof is extremely similar to the proof of 
Theorem 7.  As in that case, it is sufficient to establish (5.20a)
for an arbitrary continuous function $f$. Such a function can be 
written as an absolutely convergent sum of Fourier modes (with respect
to rotations).  The coefficient of each mode is a $\C^2$ valued
function of $s$.  The operator $e^{t{L^*}}$ commutes with rotation, and
so acts separately on each Fourier mode.

On the $n$-th Fourier mode, $e^{t{L^*}}$, for $t$ small, acts as follows.
$$ \left( e^{t{L^*}} f \right) (s) = \cases{
A \psi(s-t+\a) + B \psi(s-t+\b) & if $s \in [0,t)$; \cr
\psi(s-t) & all other $s \in [0,\mu)$,} 
\eqno(5.21) $$
where the matrices $A$ and $B$ are, as in \S 4,
$$ A = \left ( \matrix{0 & \exp[-in(\theta + \pi/2)] \cr
\exp[in(\theta + \pi/2)] & 0} \right ), \eqno(5.22)
$$
$$
B = \left ( \matrix{\exp(in\theta) + \exp(i n (\theta+\pi)) & 2 \exp(-in\theta) \cr
2 \exp(in\theta) &  \exp(-in\theta) + \exp(i n (-\theta+\pi))} \right ).
\eqno(5.23) $$

The 0-th Fourier mode decouples into $1 \choose 1$
and $1 \choose -1$ components.  The $1 \choose 1$ component is 
the distribution of sizes regardless of orientation, and its asymptotic
behavior was already computed in Theorem 8. We must show that the
$n=0$ $1 \choose -1$ component, and all the Fourier modes with $n \ne 0$,
grow strictly slower than the size of the system, and so
represent a decreasing fraction of the system. 

We will control the $L^1$ norms of the unwanted Fourier modes.  To
do this we need the $L^1$ norms of the matrices $A$ and $B$, and
various products of $A$ and $B$.  The $L^1$ norm of a matrix is 
maximum, over all columns, of the sum of the absolute values of 
the entries in that column.  One can get a bound on the growth of
the $L^1$ norm of a mode by the mode with its absolute value, and
replacing the matrices $A$ and $B$ by their norms.

For $n=0$, $B {1 \choose -1}=0$.  With $B=0$ it is as if there is only
one daughter per division, hence the $L^1$ norm of the $1 \choose -1$
mode at time $t$ is bounded by the $L^1$ norm of the mode at time 0.
Hence, as a fraction of the system, this mode shrinks like $e^{-2s}$.

Next we consider $n$ odd, for which the diagonal terms in $B$ vanish.  
The sum of the absolute values of the entries of each column of 
$B$ equals two.  
This is as if, at each subdivision, only two daughter B tiles are
produced, instead of 4.  To put it another way, at each subdivision
a fraction $b^2/2c^2$ of the area
is lost.  Since each piece of a tile of size $-t$
must be divided at least $t/\b$ times, this means that the $L^1$ norm of
the $n$-th mode, for $n$ odd, can grow no faster than
$e^{2t} [1 - (b^2/2c^2)]^{t/\b}$ and so, as a fraction of the system,
goes to zero.

Finally we consider $n$ even but nonzero.  Here the column sums of 
$A$ and $B$ are the same as in the $n=0$ case, namely 1 and 4, respectively.
However, the $L^1$ norms of various products of $A$ and $B$ are
smaller that in the $n=0$ mode.  For example, 
$ B^2 = 4 \cos(n\theta) B$ has norm $16 |\cos(n \theta)|$, which is
strictly smaller than 16.   The norm of $BAB$ is also 
$16 |\cos(n\theta)|$.  Indeed, the only
words in $A$ and $B$ which have norms as large as in the $n=0$
case are $A^m$ and $A^{m_1} B A^{m_2}$.  Since the expansion of $e^{t{L^*}}f$,
for $t>\b$, involves expressions such as $B^2$, the growth of 
the $L^1$ norm of 
$n$-th Fourier mode is bounded by an exponent strictly less than 2. 
As a fraction of the system, the $n$-th mode goes to zero.  \qed

\nd {\bf Remark:} The spectrum of $e^{t{L^*}}$ may be obtained exactly as
in \S 4.  In seeking eigenvectors, equation (4.20) is replaced by
$$ (e^{\mu\l} - e^{(\mu-\b)\l}B - e^{(\mu-\a)\l}A)v=0. \eqno(5.24) $$ 
Subsequent analysis may be repeated word for word, replacing $\l$ by
$e^\l$, $p$ by $\a$, $q$ by $\b$, $m$ by $\mu$, $\psi_k$ by $\psi(s)$,
and $r^{k}$ by $e^{-s}$.

{\nd \bf \S  6.  Two exceptional tilings --- $Til(1/3)$ and $Til(2)$}

We saw in \S 4 how the population statistics of a rational tiling
$Til(p/q)$ depends on $p$ and $q$.  If $q>1$, the second eigenvalue
of the population matrix is greater than one, and fluctuations
increase with subdivision.  This leads to phenomena such as slippage
along fault lines and a failure to be globally edge-to-edge. 
A typical example, $Til(1/2)$, was studied in \S 3.

In this section we study examples of the remaining cases.  We study
$Til(2)$ as an example of a $Til(p/1)$ tiling.  In all such tilings,
the second eigenvalue of the population matrix is less than
one.  We shall see how, in the case of $Til(2)$, this prevents
slippage along fault lines.

\vs .1
\hs1.15 \vbox{\epsfxsize=3truein\epsfbox{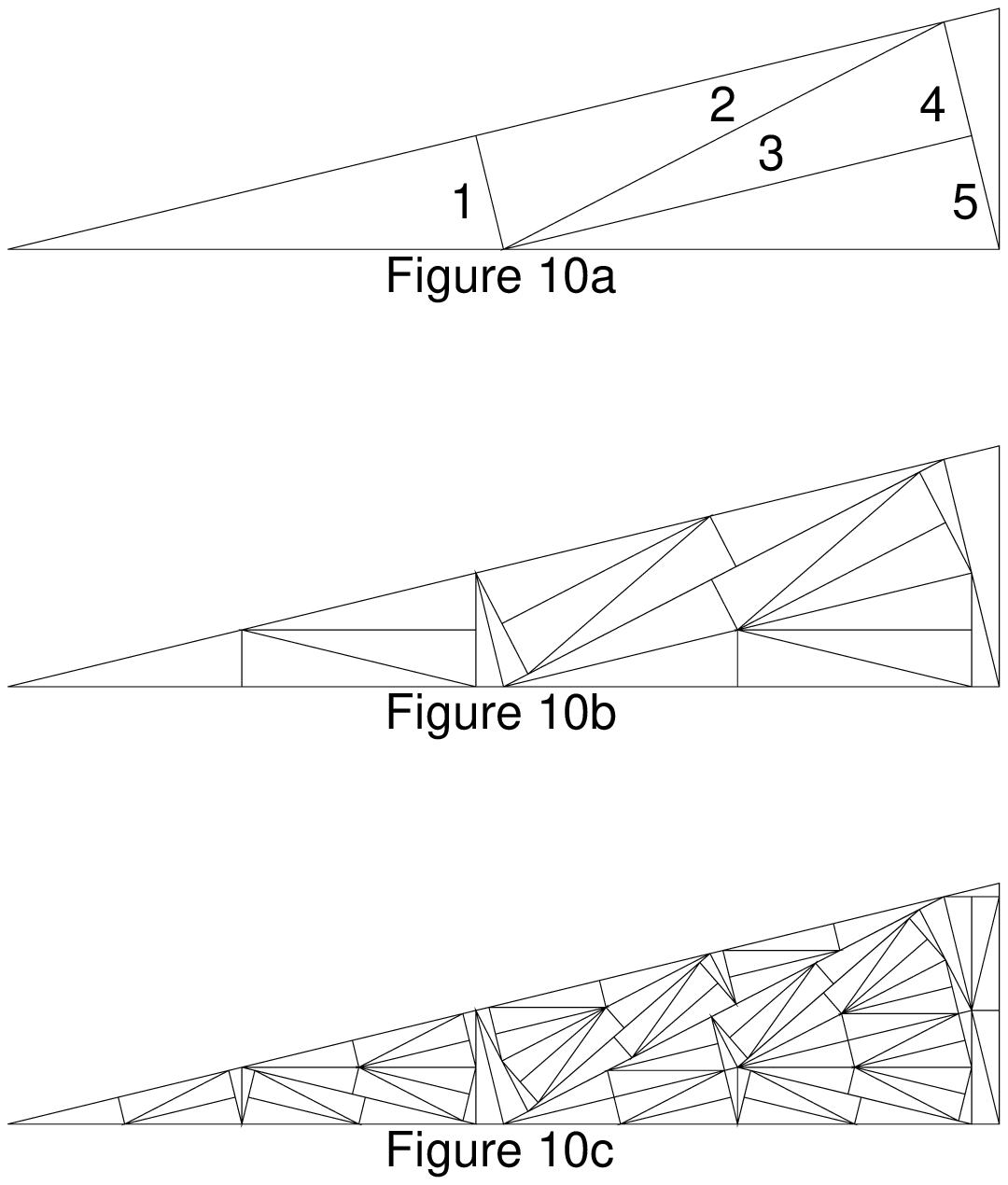}} 
\vs .1
\centerline{Figure 10. Three stages of subdivision for $Til(2)$.}
\vs .1

Finally, we consider $Til(1/3)$, the only rational 
tiling to exhibit only a finite number of orientations of each size 
of tile.  As in all cases with $q>1$, there is an eigenvalue greater
than one in the problem.  The fluctuations governed by this eigenvalue
are enough to force $Til(1/3)$ to be nonperiodic.  However, because
of rational relations between the lengths of certain edges, the tiles
in $Til(1/3)$ meet in only a finite number of ways.

We begin with $Til(2)$.  $Til(2)$ is based on the right triangle with
$a= \sqrt{5}-2 \approx 0.2361$, $b = 2 \sqrt{\sqrt{5}-2} \approx 0.9717$,
$c=1$.  Several iterations of the subdivision are shown in Figure 10.
An essential feature of $Til(2)$ is

\nd {\bf Theorem 10} {\it The triangles in the tiling $Til(2)$ meet
in only a finite number of ways.}

\nd {\it Proof:} The proof is essentially in two steps.  First we show
that the slippage along the primary fault line is bounded.  Then we show
that slippage along a fault line is the {\it only} means by which tiles in
a rational tiling can meet in an infinite number of ways, and that all
fault lines are similar to the primary fault line. 

To examine what happens along the primary fault line, we 
consider the boundary of $T_{2n}$, the $2n$-th subdivision of
the basic triangle.  Note that the hypotenuse and short
leg of $T_2$ consist only of hypotenuses and short legs of big triangles.
Applying the subdivision again, we get that the hypotenuse and short
leg of $T_4$ also consists only of hypotenuses and short legs of big 
triangles.  Similarly for all $T_{2n}$.  

As in the proof of Theorem 5, we consider the evolution of these legs
and hypotenuses as a one dimensional 
substitution system in its own right.  Let $\sigma$ denote the
effect of subdividing the basic triangle twice.
Under $\sigma$, each hypotenuse (denoted $H$) 
is replaced by 4 hypotenuses and
a short leg (denoted $S$), while each short leg is replaced by a hypotenuse.
That is, the one dimensional population matrix is
$$ M= \left (\matrix{4 & 1 \cr 1 & 0} \right ), \eqno (6.1) $$
with eigenvalues $\lambda_\pm = 2 \pm \sqrt{5}$ and eigenvectors
$v_\pm = \left ( \matrix{\lambda_\pm \cr 1} \right )$. 
 
Now let $H_n$ and $S_n$ be the number of hypotenuses and short legs in
$\sigma^n(H)$, and let $H_n'$ and $S_n'$ be the number of hypotenuses
and short legs in $\sigma^n(S)$.  By expanding $0 \choose 1$ and $1
\choose 0$ in terms of $v_\pm$, it is easy to see that
$$ (\sqrt{5}-2) H_n - S_n = -(2-\sqrt{5})^{n+1}; \qquad
 (\sqrt{5}-2)H_n' -  S_n' = -(2-\sqrt{5})^{n}. \eqno(6.2) 
$$

Next we measure slippage.  Let $P$ and $R$ be the endpoints of a
hypotenuse, as in Figure 7, and let $E$ be any intermediate point, not
necessarily the midpoint.  Let $f_n(E)$ be the number of complete
short legs, between $P$ and $E$, in $\sigma^n(PR)$, minus the number
of short legs between $P$ and $E$ in $\sigma^n(RP)$. As in the proof
of Theorem 5, $f_n(E)$ measures the extent to which the two tiles of
$T_{2n+2}$ that meet at $E$ are offset.

\nd {\bf Lemma 10} {\it $|f_n(E)| < 5$}. 

\nd {\it Proof:} Let $v_0=P$, and let $v_k$ be the vertex of 
$\sigma^k(PR)$, between $P$ and $E$, that is closest to $E$.  Note
that, in $\sigma^k(PR)$, there are at most 4 hypotenuses and at most
one short leg between $v_{k-1}$ and $v_k$, since the interval
$v_{k-1}v_k$ was only part of a hypotenuse or short leg in
$\sigma^{k-1}(PR)$.  By (6.2), $(\sqrt{5}-2)$ times the number of
hypotenuses in $\sigma^n(PR)$ between $v_{k-1}$ and $v_k$, minus the
number of short legs, is bounded in absolute value by
$(\sqrt{5}-2)^{n-k}$.  Summing over $k$, we get that $(\sqrt{5}-2)$
times the number of complete hypotenuses between $P$ and $E$ in
$\sigma^n(PR)$, minus the number of short legs, is bounded in absolute
value by $\sum_{i= 0}^\infty (\sqrt{5}-2)^i < 2$.  A similar bound
applies to the number of hypotenuses and legs in $\sigma^n(RP)$.  Thus
the surplus of short legs on one side of $PE$ relative to the other,
plus $(\sqrt{5}-2)$ times the deficit of hypotenuses, is bounded by
$2+2=4$.  Since a surplus of short legs implies a deficit of
hypotenuses, the surplus of short legs is itself bounded by 4.  \qed

We return to the proof of Theorem 10.  Lemma 10 limits the number of
ways for two tiles to meet across a fault line.  Suppose two tiles
$t_1$ and $t_2$ meet across a fault line $PR$, modeled on
$\sigma^n(H)$.  Let $E$ be a point on their common edge.  The distance
from the vertex of $t_1$ closest to $P$ to the vertex of $t_2$ closest
to $P$ is either $|f_n(E)|(\sqrt{5}-2)$ or $1-|f_n(E)|(\sqrt{5}-2)$.
Since $|f_n|$ is at most 4, this means there are only a finite number
of ways for two triangles to meet across such fault lines.  Now
$\sigma^n(H)$ is the result of subdividing the basic triangle an even
number of times.  However, since subdivision is deterministic, having
only a finite number of distinct configurations in the even
subdivisions implies that there are only a finite number of distinct
configurations in the odd subdivisions, and thus a finite number of
configurations in all.

To complete the proof of Theorem 10, we must show that every pair of
adjacent triangles either meets full-face to full-face, or meets
across a fault line based on successive subdivision of a hypotenuse.
Consider two tiles, $t_1$ and $t_2$, that meet. Let $k$ be the
smallest integer such that both tiles lie in the same supertile $S$ of
order $k$.  Since $t_1$ and $t_2$ do not meet in a supertile of order
$k-1$, $t_1$ and $t_2$ must meet across one of the five lines of the
first subdivision of $S$.  See Figure 10a.  There is a local
reflection symmetry across edges 1 and 3, so if $t_1$ and $t_2$ meet
across these edges they must meed full-face to full-face.  Edge 2 is
the primary fault line.  Further division (see Figure 10c) shows that
there is local reflection symmetry across edge 4, while edge 5 is a
hypotenuse-based fault line, as considered above.  \qed

We now turn to the tiling $Til(1/3)$.  By Theorem 4, $Til(1/3)$ is the
only tiling in our construction to have both a finite number of sizes
of tiles, each of which appears only in a finite number of
orientations.  $Til(1/3)$ is based on an isosceles right triangle.
Although the two legs have the same length, we distinguish between the
two, calling the ``$b$'' side ``long'' and the ``$a$'' side ``short'',
in analogy to the tilings with $b>a$.  In subdividing we must specify
which legs of the daughter tiles are labeled ``long'' and ``short''.
This is shown in Figure 11, and several further subdivisions are shown
in Figure 12. 

\vs .1
\hs1.4 \vbox{\epsfxsize=2.5truein\epsfbox{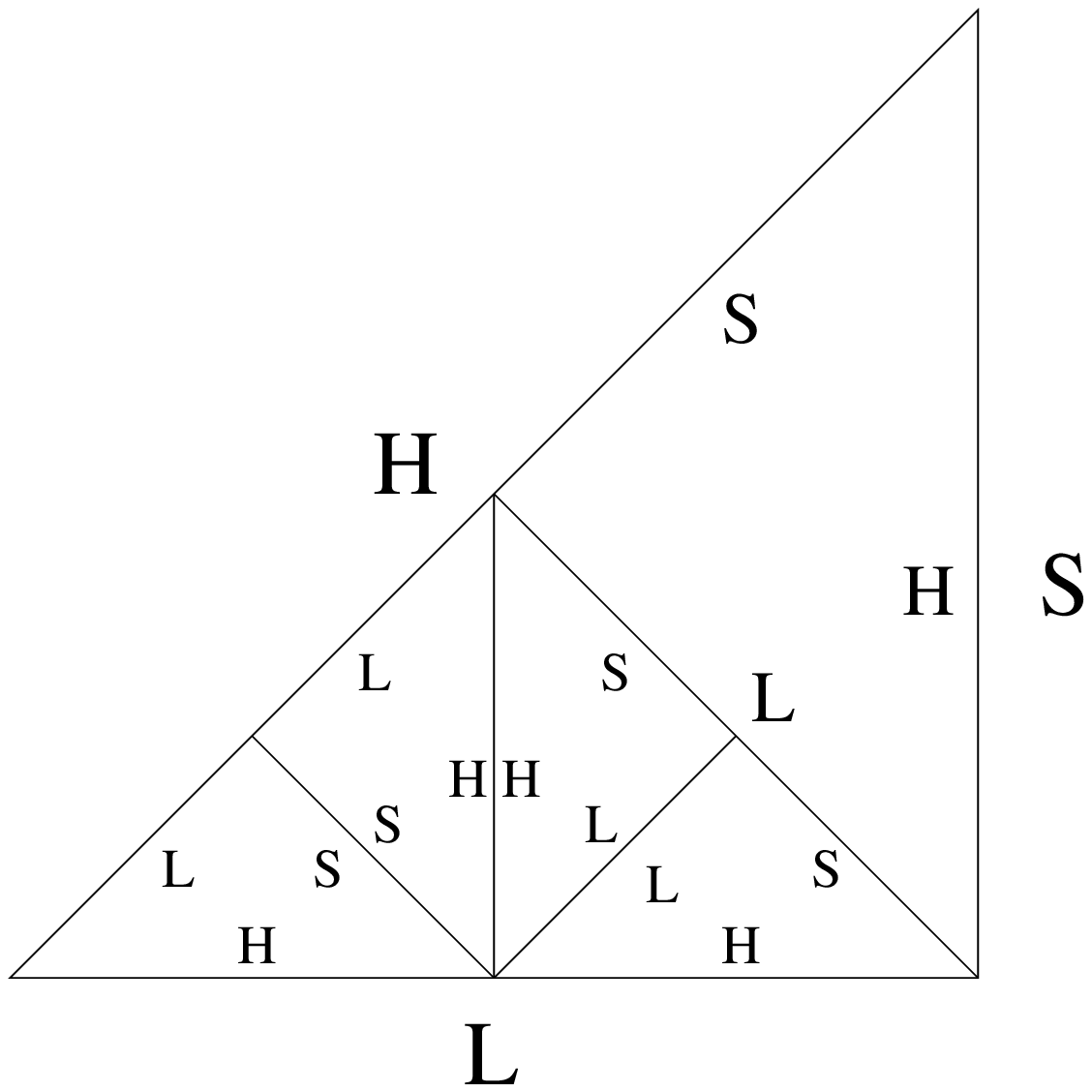}} 
\vs .1
\centerline{Figure 11. ``Long'' and ``short'' edges in $Til(1/3)$.}
\vs .1

\vs .1
\hs0.65 \vbox{\epsfxsize=4truein\epsfbox{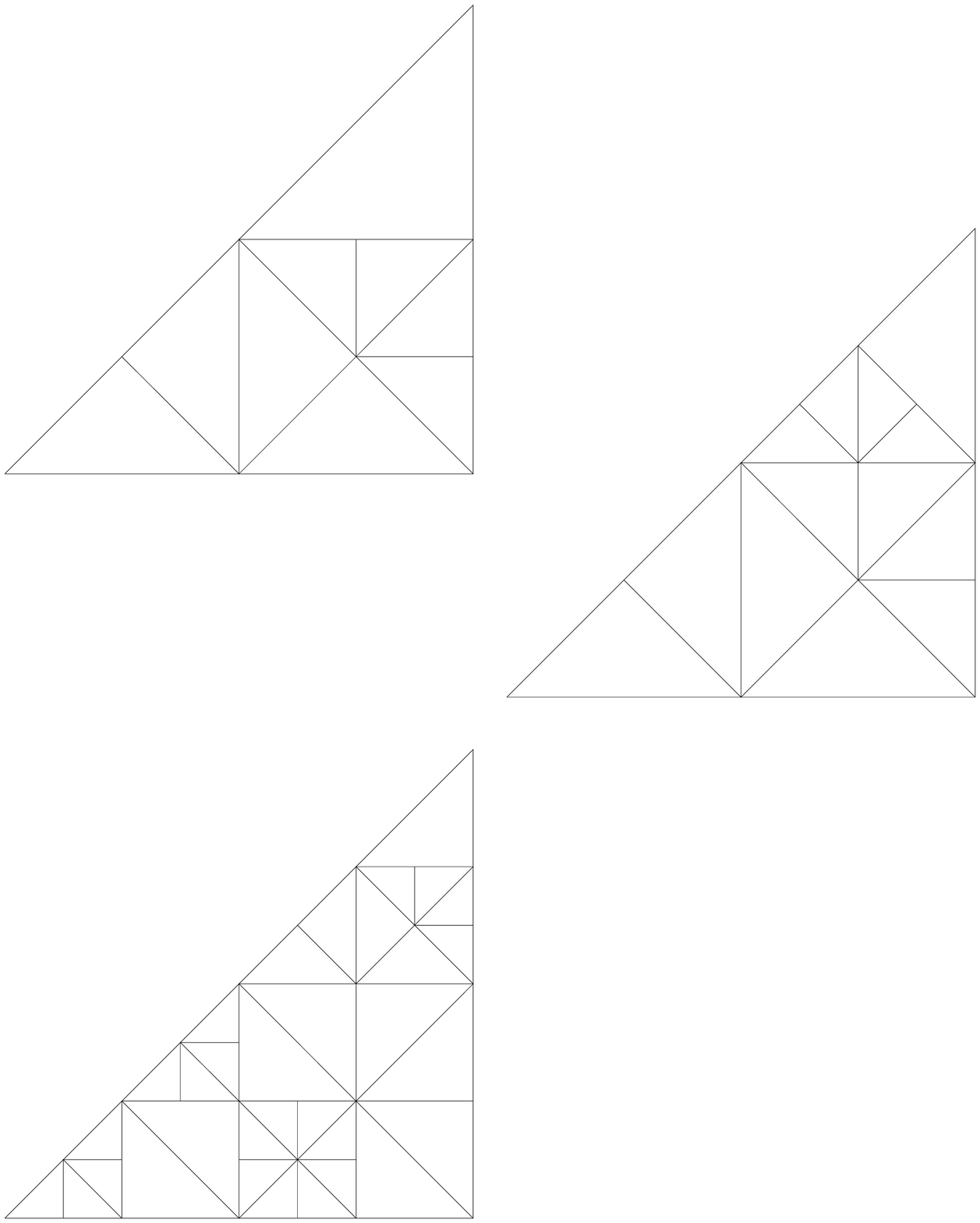}} 
\vs .1
\centerline{Figure 12. Three stages of subdivision for $Til(1/3)$.}
\vs .1

\nd {\bf Theorem 11} {\it The tiling $Til(1/3)$ is nonperiodic.  The tiles
meet in only a finite number of ways.}

\nd {\it Proof:}  As usual, we consider the 1-dimensional substitution
scheme induced on the edges by 2-fold substitution.  The $2n$-th
subdivision of the hypotenuse consists of some hypotenuses of large
triangles (denoted $H$), some long legs of medium sized triangles
($L$), and some hypotenuses of small triangles ($h$).  Let $\sigma$
denote the action of two subdivisions.  $\sigma$ takes each $H$ to an
$H$ and two $L$'s, each $L$ to two $h$'s, and each $h$ to an $H$.  The
population matrix is
$$ M = \left ( \matrix{ 1 & 0 & 1 \cr 2 & 0 & 0 \cr 0 & 2 & 0} \right )
\eqno (6.3) 
$$
with eigenvalues $2$ and $(-1 \pm \sqrt{-7})/2$, and with the eigenvector
$\left( \matrix{1 \cr 1 \cr 1} \right )$ corresponding to eigenvalue 2.
Note that the complex eigenvalues have magnitude $\sqrt{2}$.

Asymptotically, the three types of edges appear in a ratio of 1:1:1.
However, since the initial condition $\left( \matrix{1 \cr 0 \cr 0}
\right )$ is not an eigenvector, and since all eigenvalues are greater
than 1 (in magnitude), the difference in number between $H$'s and
$L$'s, or $L$'s and $h$'s, with grow exponentially with time.  That
is, while the total population grows as $2^n$, the fluctuations grow
as $\sqrt{2}^n$.

If $Til(1/3)$ were periodic, a long
line of the form $\sigma^n(H)$ would consist of several periods, plus
a remainder at each end. Each period would have $H$'s, $L$'s, and
$h$'s in exactly a 1:1:1 ratio, so only the partial periods at each
end could contribute to the difference in population between $H$ and $L$. 
Thus the population difference would remain bounded as $n \to \infty$.
Since this difference is unbounded, $Til(1/3)$ cannot be periodic.

Finally we note that the length of $h$ is the same as that of $L$, and
half that of $H$.  This simple ratio of lengths means that slippage
along the fault line has no effect on the number of ways triangles can
meet.  If two tiles meet across $\sigma^n(H)$, either they have a
vertex in common or their closest vertices are separated by the length
of $h$. Thus there are only a finite number of ways for two tiles to
meet across a fault line.  By the same argument as in the proof of
Theorem 10, this implies that there are only a finite number of ways
for triangles to meet at all. \qed
 
\vs .20

{\nd \bf \S  7.  Conclusions}

We have constructed a family of substitution systems, indexed by the
parameter $z = \ln(\sin(\theta))/\ln[\cos(\theta)/2]$, where $\theta$
is an angle in the basic triangular tile.  We have established the
following properties.

\item {1)} The tilings generated by these substitutions are all 
non-periodic.
\item {2)} The tilings have well-defined limiting distributions of
size and orientation.  If $\theta/\pi$ is irrational, this
distribution is rotationally invariant.  In Radin's terminology, the
tilings have ``statistical rotational symmetry''.  The form of the
joint distribution of size and orientation suggests that the tiling
has a purely absolutely continuous spectrum.
\item {3)} The tilings with rational $z$ all satisfy the hypotheses of 
Goodman-Strauss's theorem, implying that they can be forced through
local matching rules.
\item {4)} The rational tilings $Til(p/q)$, with $q>1$, have statistical
fluctuations that grow with iterations of the substitution rule (although
they grow slower than the size of the system).  
In $Til(1/2)$ these fluctuations force triangles to meet
in an infinite number of distinct ways.  We conjecture that this infinite
diversity of local behavior is a property of all rational tilings $Til(p/q)$ 
with $q>1$ and $p/q \ne 1/3$.
\item {5)} In the rational tilings $Til(p/1)$ the eigenvalues that control
fluctuations are all less than one.  In $Til(2)$ this forces the tiles
to meet in only a finite number of local patterns.  We conjecture
that this is a property of all tilings $Til(p/1)$.

\vs.20 \nd
{\bf Acknowledgements.}\  It is a pleasure to thank Chaim Goodman-Strauss,
Yoram Last, Tom Mrowka, Johan R\aa{}de, Charles Radin and 
Felipe Voloch for useful discussions, and Steven Janowsky for assistance
with computer graphics.

\vs.3

\nd {\bf References}
\vs.3 \nd
[CFS]\ I. Cornfeld, S. Fomin and Ya. Sinai, ``Ergodic Theory'',
Springer-Verlag, New York, 1992
\vs.1 \nd
[GS]\ C. Goodman-Strauss, {\it Matching rules and substitution tilings},
preprint 1996.
\vs.1 \nd
[R1]\ C.\ Radin, {\it The pinwheel tilings of the plane}, {Annals of
Math.}  {\bf 139} (1994), 661--702.
\vs.1 \nd
[R2]\ C.\ Radin, {\it Space Tilings and Substitutions}, {Geometriae Dedicata}  {\bf 55} (1995), 257--264.
\vs.1 \nd
[R3]\ C.\ Radin, {\it Symmetry and Tilings}, {Notices of the AMS}  {\bf 42} (1995), 26--31.
\vs.1 \nd
[Sen]\ M. Senechal, ``Quasicrystals and geometry'', Cambridge
University Press, Cambridge, 1995
\vfill
\end